\documentclass[fleqn,preprint,authoryear,5p,serif]{elsarticle}
\usepackage{amsmath,amssymb}
\usepackage{mathrsfs}
\usepackage{color}
\usepackage{caption}
\usepackage{booktabs}
\usepackage{multirow,booktabs}
\usepackage{subfigure}
\journal{Journal publication}
\usepackage[normalem]{ulem}
\usepackage{enumitem}
\newlist{steps}{enumerate}{1}
\setlist[steps, 1]{label = Step \arabic*:}

\def\ee		{\mathrm{e}}
\def\ii		{\mathrm{i}}
\def\diff	{\mathrm{d}}

\def\Re		{\mathrm{Re}}
\def\Im		{\mathrm{Im}}

\def\XXint#1#2#3{{
		\setbox0=\hbox{$#1{#2#3}{\int}$}
		\vcenter{\hbox{$#2#3$}}\kern-.5\wd0}}

\newcommand{\varendash}[1][5pt]{%
	\makebox[#1]{\leaders\hbox{--}\hfill\kern0pt}%
}

\newlength{\intwidth}

\begin{document}
\begin{frontmatter}
\title{Accurate and efficient hydrodynamic analysis of structures with sharp edges by the Extended Finite Element Method (XFEM): 2D studies}

\author[HEU]{Ying Wang}

\author[DTU]{Yanlin Shao\corref{cor}}
\cortext[cor]{Corresponding author}\ead{yshao@mek.dtu.dk}

\author[HEU]{Jikang Chen}

\author[TCOMS]{Hui Liang}

\address[HEU]
{College of Shipbuilding Engineering, Harbin Engineering University, Harbin, 150001, China}

\address[DTU]
{Department of Mechanical Engineering, Technical University of Denmark, 2800 Lyngby, Denmark}

\address[TCOMS]
{Technology Centre for Offshore and Marine, Singapore (TCOMS), 118411, Singapore}

\begin{abstract}
	
  \par Achieving accurate numerical results of hydrodynamic loads based on the potential-flow theory is very challenging for structures with sharp edges, due to the singular behavior of the local-flow velocities. In this paper, we introduce, perhaps the first time in the literature on marine hydrodynamics, the Extended Finite Element Method (XFEM) to solve fluid-structure interaction problems involving sharp edges on structures. Compared with the conventional FEMs, the singular basis functions are introduced in XFEM through the local construction of shape functions of the finite elements. Four different FEM solvers, including conventional linear and quadratic FEMs as well as their corresponding XFEM versions with local enrichment by singular basis functions at sharp edges, are implemented and compared. To demonstrate the accuracy and efficiency of the XFEMs, a thin flat plate in an infinite fluid domain and a forced heaving rectangle at the free surface, both in two dimensions, will be studied. For the flat plate, the mesh convergence studies are carried out for both the velocity potential in the fluid domain and the added mass, and the XFEMs show apparent advantages thanks to their local enhancement at the sharp edges. Three different enrichment strategies are also compared, and suggestions will be made for the practical implementation of the XFEM. For the forced heaving rectangle, the linear and 2nd order mean wave loads are studied. Our results confirm the previous conclusion in the literature that it is not difficult for a conventional numerical model to obtain convergent results for added mass and damping coefficients. However, when the 2nd order mean wave loads requiring the computation of velocity components are calculated via direct pressure integration, the influence of singularity is significant, and it takes a tremendously large number of elements for the conventional FEMs to get convergent results. On the contrary, the numerical results of XFEMs converge rapidly even with very coarse meshes, especially for the quadratic XFEM. Unlike other methods based on domain decomposition when dealing with singularities, the FEM framework is more flexible to include the singular functions in local approximations. 
\end{abstract}


\end{frontmatter}

\section{Introduction}

  \par Numerical analysis is playing an increasingly important role in marine hydrodynamics. Computational Fluid Dynamic (CFD) models based on the Navier-Storkes (NS) equations with proper turbulence modeling are the most comprehensive ones for this purpose. They are applicable in more applications than a potential-flow model, in particular when viscous flow separation and wave breaking become relevant and important. The computational costs, however, are normally too high to afford, which is regarded as one of the bottlenecks of CFD models, if they are heavily involved in the design of marine structures. 
  Due to large-volume nature of most of the marine structures, the inertial effect is predominant whereas viscosity effect plays a secondary role. Therefore, the potential-flow theory is often applied together with empirical corrections for viscous effects.
  
  \par For the potential-flow problems, Boundary Element Method (BEM) is the most commonly used numerical method in marine hydrodynamics, as it can reduce the dimension of the problem by one and only the boundaries of the fluid domain need to be discretized. Even though the number of unknowns is reduced in BEM compared with a volume method, it is still challenging for a conventional BEM to solve the resulting linear system with a large number of unknowns, because the matrix is dense.  $O(N^2)$ memory is required by the conventional BEMs, and $O(N^2)$ and $O(N^3)$ operations are required for iterative solvers and direct solvers, respectively. Here $N$ denotes the number of total unknowns on the boundary surfaces. 
  
  \par Although BEM is a very popular numerical method in potential-flow hydrodynamic analyses, field solvers are also widely used. \cite{wu1994finite} is among the first to use FEM to investigate 2D nonlinear free-surface flow problems in the time domain. \cite{wu1995time} studied the fully-nonlinear wave-making problem by both FEM and BEM, and suggested that FEM is more efficient than BEM in terms of both CPU time and computer memory. \cite{2010FinitePart1, 2010FinitePart2} used a FEM to simulate the interaction between 3D fixed bodies and steep waves. On the other hand, high-order volume methods have gained great interest. \cite{2007On} and \cite{engsig2009efficient} developed 2D and 3D high-order Finite Difference Methods (FDMs) to study fully-nonlinear water wave problems in potential flows. \cite{2012Shao} and \cite{shao2014harmonic} proposed high-order Harmonic Polynomial Cell (HPC) methods in 2D and 3D respectively to study water waves and their interaction with structures. Some recent extensions were made to utilize immersed boundary strategies and overset meshes to achieve better accuracy and stability \citep[e.g., see][]{hanssen2018free, tong2019numerical, tong2021adaptive,law2020numerical,liang2020liquid}. Compared to the BEMs, field solvers deal with sparse matrices, and the computational costs are roughly linearly dependent on the number of unknowns. 
  
 Ordinary boundary-element and volume methods, e.g. BEM, FEM, FDM and HPC methods, are based on local approximations using smooth functions. Thus, very fine meshes have to be applied at areas where the fluid solution tends to be singular. Sharp edges are widely present in typical offshore structures. Examples are pontoons of semi-submersibles and tension leg platforms \citep[e.g., see][]{chen1995numerical, zhou2015resonance}, damping plates on offshore platforms \citep[e.g., see][]{tao2007spacing, shao2016stochastic, shao2019pontoon} and offshore floating wind turbine structures \citep[e.g., see][]{xu2019effect}, as well as the bilge keels on the ships. Besides, the analytical methods, such as the multi-term Galerkin method \citep[e.g., see][]{2019New,1995complementary}, have also been used to include the local singularities. From industrial application point of view, it is essential to be able to obtain accurate numerical results with affordable computational efforts. However, this is not always possible, even for the 2nd order mean wave loads. 
  
 The calculation of 2nd order mean wave loads involves quadratic terms of the 1st order quantities, which pose great challenges at the sharp edges where the fluid velocities tend to be infinite. \cite{taylor1993effect} investigated the effect of corners on diffraction/radiation wave loads and wave-drift damping, and revealed that the most important hydrodynamic loads and the amplitudes of body motion do not change significantly while the radius of the corner approaches zero. For a floating truncated vertical cylinder free to surge and heave, \cite{zhao1989interaction} found it is difficult to obtain convergent 2nd order mean wave forces via the direct pressure integration. In their work, a method based on momentum and energy relationship was shown to be more robust and efficient. By applying the variants of Stokes's theorem, \cite{dai2005computation} and \cite{chen2007middle} developed a `middle-field formulation', which transforms the body-surface integral to a control surface at a distance from the body. Similar strategy was also applied by  \cite{liang2017multi-domain} where a multi-domain approach was developed. The middle-field formulation can be used to calculate drift forces and moments in all 6 degrees of freedom. The floating truncated vertical cylinder studied by \cite{zhao1989interaction} was revisited in \cite{2018Numerical} and four different methods were used to calculate the vertical mean wave force, including a momentum formulation implemented in a time-domain higher-order BEM \citep{shao2013second}, a semi-analytical solution \citep{mavrakos1988}, the middle-field method in HydroStar, as well as the near-field method in HydroStar. The first three methods matched very well with each other, confirming the accuracy of the earlier results by \cite{zhao1989interaction} based on momentum and energy relations. However, the results determined by the direct pressure integration were quite different in the heave resonance regime. As elucidated in \cite{2018Numerical}, the results by the direct pressure integration are not convergent, despite very fine meshes have been used. 
 
 \cite{2020Comparative} used five different methods to investigate nonlinear radiation forces of bodies with sharp or rounded edges in the time domain. The first four methods are all near-field methods, and the fifth one based on momentum conservation. They found that the singularity at the sharp edge plays significant roles on numerical computation of hydrodynamic forces in all near-field methods, while it has much less influences on results based on momentum conservation. Using an approach based on a control surface, \cite{2020A} rewrote the integration of velocity square terms on body surface into the sum of two other integrals, one on a control surface enclosing the structure and the other on the free surface between the structure and the control surface. Encouraging results were obtained for double-frequency wave-radiation forces on an oscillating truncated vertical cylinder.

 \par This paper aims to introduce, verify and demonstrate the XFEM as an accurate and efficient tool to calculate the linear and 2nd order wave loads on structures with sharp edges, without having to use a control surface. The XFEM has a powerful framework, which allows for adding the knowledge of the local solutions, normally known as a priori, to the finite-element approximation space at specific nodes. The solution enrichment at those nodes does not require any modification to the meshes. The idea of XFEM was originally used by \cite{1999Elastic} to solve the problem of elastic crack growth, and one year later, \cite{2000Arbitrary} formally named this approach as XFEM. The XFEM can be seen as an extension of the standard FEM based on the conception of Partition of Unity (PU) \citep{Babu1997THE}, and thus it maintains all advantages of the standard FEM. Earlier concepts of PU dates back to 1994, when it was first used to solve the so-called roughness coefficient elliptic boundary value problem by \cite{1994Special} with the name of special finite element method, namely the Generalized Finite Element Method (GFEM). Based on the ideas in \cite{1994Special}, the GFEM was further elaborated by \cite{1995On}, \cite{1996PUFEM}, \cite{I1997THE} and \cite{1997Approximation} with the name of partition of unity method (PUM) or partition of unity finite element method (PUFEM). The GFEM was developed in \cite{2000Thedesign} and \cite{2000Thegeneralized} with the name of GFEM. In the early days, both XFEM and GEFM were developed independently even their basic idea is similar. A feature to distinguish the XFEM and the GFEM in early work is that only local parts of the domain are enriched by XFEM, but GFEM enriches the whole domain globally. However, \cite{2010Fries} argued that the XFEM and the GFEM are almost identical numerical methods. The XFEM represents the singular properties by adding singular basis function or any analytical recognition of the solution to local approximation space, and it has been a tremendous success in dealing with singular or discontinuous problems, no matter how strong the discontinuity is \citep[see, e.g.][]{2001Modeling, 2002Non, 2015Extended, 2010Fries}. Besides, XFEM has also been introduced to CFD to model two-phase flows \citep{2010The}. 
  
  \par   In the present work, as verification and demonstration, flow around an infinite-thin flat plate and a heaving rectangle on the free surface will be studied via four different FEMs, namely the linear FEM, linear XFEM, quadratic FEM and quadratic XFEM. Convergence studies will be presented to illustrate accuracy and efficiency of the XFEMs. Our results indicate that the singularities at sharp edges do not have a strong influence on the calculating of added mass and damping, confirming the conclusion from an earlier study by \cite{taylor1993effect}. However, if the near-field method is used, it is extremely challenging for conventional FEMs to achieve convergent 2nd order vertical mean forces for the heaving rectangle with affordable computational time on a normal PC. On the contrary, the XFEMs with local enrichment, using corner-flow solutions \citep{Newman2017Marine} at sharp edges, can achieve convergent results with much coarser meshes. Three different local enrichment strategies of XFEM will also be compared and suggestions will be made for practical implementation.
  
  \par The rest of the paper will be organized as follows. In Sect.~\ref{sect:mathematical-formu}, the formulation of the boundary-value problem and corner-flow solutions are presented. In Sect.~\ref{sect:Numerical Method}, the basic idea of conventional FEM and the XFEM are introduced via a mixed boundary-value problem in 2D. Besides, three enrichment strategies for the XFEM are presented and compared. In Sect.~\ref{sect:Numerical case}, as the first verification case, the velocity potential in fluid domain and the added mass of an infinitely-thin flat plate are studied and compared with the analytical solution. The second verification concerns a heaving rectangle on a free surface, solved in the frequency domain. In Sect.~\ref{sect:Conclusion and perspective}, some conclusions are drawn.
  
\section{Mathematical formulation}\label{sect:mathematical-formu}

  \subsection{Governing equation and linearized boundary condition}

  \begin{figure*}[t]
  	\centering
  	\includegraphics[scale=.55]{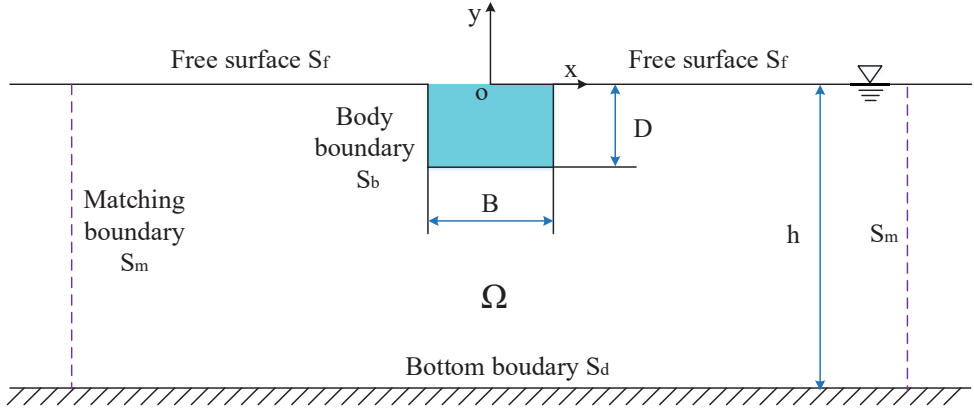}
  	\caption{An illustration of the fluid domain and its boundaries, as well as the definition of the coordinate system.}
  	\label{Fig.21}	
  \end{figure*}
  \par A 2D coordinate system $Oxy$ is defined with the $Ox$ axis coinciding with the undisturbed free surface and $Oy$ axis orienting positively upward, as illustrated in Fig.~\ref{Fig.21}. The fluid domain $\Omega$ is enclosed by the body surface $S_b$, free surface $S_f$, bottom surface $S_d$, and vertical control surfaces $S_m$ at a distance from the body. 

  \par It is assumed that the fluid is inviscid and incompressible, and the flow is irrotational so that a velocity potential $\phi$ exists. In this study, we only consider 2D flows, and thus the governing equation in the fluid domain $\Omega$ is written as
  \begin{equation}
  \label{Eq.1}
    \frac{\partial^2 \phi}{\partial x^2} 
    + \frac{\partial^2 \phi}{\partial y^2} = 0,
  \end{equation}
  where $\phi$ denotes velocity potential. Here only radiation problem is considered, and thus the impenetrable condition on the body surface is written as:
  \begin{equation}
  \label{Eq.2}
    \frac{\partial \phi}{\partial n} = {\boldsymbol v} \cdot {\boldsymbol n} \quad \text{at}\quad S_b,
  \end{equation}
  where ${\boldsymbol v}$ is the velocity of the body and ${\boldsymbol n}$ is the vector normal to the body surface pointing out of the fluid domain. Besides, the combined linearized free-surface condition is written as 
  \begin{equation}
  \label{Eq.5}
    \frac{{{\partial }^{2}}\phi }{\partial {{t}^{2}}}
    +g\frac{\partial \phi }{\partial y}=0\quad \text { at }\quad S_f.
  \end{equation}
  The bottom condition is 
  \begin{equation}
  \label{Eq.6}
    \frac{\partial\phi}{\partial n}= 0\quad \text { at }\quad S_d.
  \end{equation}

  \subsection{Linearized frequency-domain analysis}
  
  \par Assuming that the problem is time-harmonic and a steady state is reached. Therefore, velocity potential can be separated into a spatial part and temporal part as follows:
  \begin{equation}
  \label{Eq.35}
  	\phi(x,y,t) = \Re\{\varphi(x,y) \cdot \ee^{\ii\omega t}\},
  \end{equation}
  where $\omega$ denotes the angular frequency of oscillation, and $\ii =\sqrt{-1}$. The motion of body in $j$-th mode can be defined  as:
  \begin{equation}
  \label{Eq.53}
    {\eta}_j = \Re\{\eta_{ja} \ee^{\ii\omega t}\} \quad (j=1,2,3),
  \end{equation}
  where $\eta_{ja}$ denote the amplitude of body in $j$-th mode, and $j=1$, $2$, and $3$ correspond to sway, heave, and roll motions, respectively.
  Accordingly, the governing equation and boundary-value problem (BVP) with respect to the complex velocity potential $\varphi(x,y)$ can be written as:
  \begin{equation}
  \label{Eq.36}
   \begin{aligned}
     & \frac{\partial^{2} \varphi}{\partial x^{2}}+\frac{\partial^{2} \varphi}{\partial y^{2}}=0 &\text{in}\quad \Omega, \\
     & -\omega^2 \varphi+ g\frac{\partial \varphi}{\partial y}=0 &\text{at}\quad S_f, \\
     & \frac{\partial \varphi }{\partial n} =\sum\limits_{j=1}^{3}{\ii \omega{{\eta }_{ja}}{{n}_{j}}} &\text{at}\quad S_b,\\
     & \frac{\partial \varphi}{\partial y}=0 &\text {at} \quad S_d.
  \end{aligned}
  \end{equation}
  Here $n_j$ represent the component of the normal vector in the direction of the motion of body in $j$-th mode. The dispersion relation in finite water depth is $k\tanh kh={\omega }^{2}/{g}$, where $k$ is wavenumber. Thus, the free-surface condition in Eq.~\eqref{Eq.36} can be rewritten as 
  \begin{equation}
  \label{Eq.37}
  	-(k\tanh kh) \cdot \varphi+ \frac{\partial \varphi}{\partial y}=0 \quad \text { at } y=0.
  \end{equation}
  Besides, radiation condition requiring radiated waves propagating outwards can be expressed as:
  \begin{equation}
  \label{Eq.38}
  \begin{aligned}
    \frac{\partial \varphi}{\partial x}+\mathrm{i} k\, \mathrm{sgn}(x) \,\varphi \rightarrow 0 \quad \text { when } \quad x \rightarrow+\infty.
  \end{aligned}
  \end{equation}
  If the horizontal distance between rectangle and matching boundary is large enough, the radiation condition can be satisfied at the matching boundaries $S_m$:
  \begin{equation}
  \label{Eq.39}
  \begin{aligned}
    &\frac{\partial \varphi}{\partial x}+\mathrm{i} k\, \mathrm{sgn}(x)\, \varphi = 0 \quad \text { at }  S_m.
  \end{aligned}
  \end{equation}
  
  \subsection{Corner-flow solution}\label{sect:corner-flow}
  \begin{figure}[t]
  	\centering
  	\includegraphics[scale=.40]{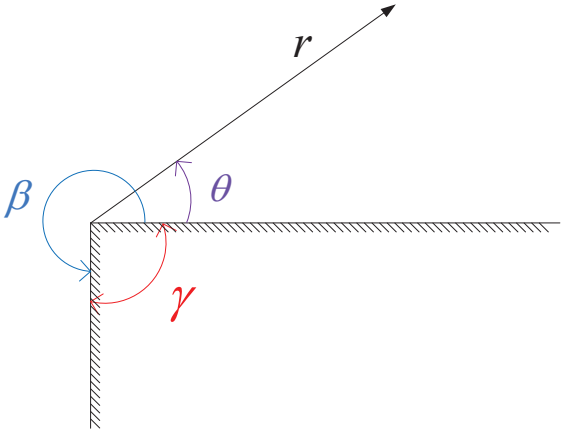}
  	\caption{Definition of the Cartesian and polar coordinate systems for the corner flow problem.}
  	\label{Fig.1}	
  \end{figure}
  \par In order to demonstrate the singular characteristics of the corner flow by potential-flow theory, the flow past a sharp corner with an exterior angle $\beta$ and the corresponding interior angle of $\gamma = 2\pi-\beta$ as shown in Fig.~\ref{Fig.1} is considered. If the considered semi-infinite wedge is fixed, the corner-flow solution can, according to \cite{Newman2017Marine}, be defined in the polar coordinate system $Or\theta$ as
  \begin{equation}
   \begin{aligned}
  \label{Eq.7}
    \varphi &=\sum\limits_{j}{{{A}_{j}}{{r}^{j\pi /\beta }}\cos \left( \frac{j\pi }{\beta }\theta  \right)}\\
    &=\sum\limits_{j}{{{A}_{j}}{{r}^{j\pi /(2\pi -\gamma )}}\cos \left( \frac{j\pi }{2\pi -\gamma }\theta  \right)},
  \end{aligned}
  \end{equation}
  where $A_j$ is a constant and $j$ is an non-negative integer number. It is obvious that the velocity determined by Eq.~\eqref{Eq.7} is singular at the tip of the semi-infinite wedge when $j\ge1$ and $\gamma<\pi$. If we define 
  \begin{equation}
  \label{Eq.8}
  m_j = \frac{j\pi }{2\pi -\gamma },
  \end{equation}
  Eq.~\eqref{Eq.7} can be rewritten as
  \begin{align}
  \label{Eq.9}
   \varphi =\sum\limits_{j}{{{A}_{j}}{{r}^{m_j}}\cos \left( m_j\theta  \right)}.
  \end{align}
  \par For a general 2D radiation-diffraction problem, the local scatter velocity potential (incident wave potential excluded) close to  a sharp edge can be expressed as
  \begin{equation}
  \label{Eq.80}
      \varphi_s =\varphi_0+\sum\limits_{j}{{{A}_{j}}{{r}^{m_j}}\cos \left( m_j\theta  \right)}.
  \end{equation}
  Here the first term $\varphi_0$ is a smooth velocity potential satisfying the non-trivial Neumann-boundary conditions for the scatter velocity potential 
  \begin{equation}
  \label{Eq.81}
      \frac{\partial \varphi_s}{\partial n} = -\frac{\partial \varphi_i}{\partial n} + \mathbf{v} \cdot \mathbf{n},
  \end{equation}
  where $\varphi_i$ is the incident wave potential, $\mathbf{v}$ is the rigid-body velocity at the edge and $\mathbf{n}$ is the normal vector on the body surface. For a radiation problem, an example of $\varphi_0$ has been given by \cite{liang2015application} as $\varphi_0 = u\,x + v\,y$. $u$ and $v$ are the horizontal and vertical velocities at the corner due to rigid-body motions, respectively. The second term of Eq.~\eqref{Eq.80} represents the corner-flow solutions derived from zero Neumann-boundary condition at the sharp edges. Since the first term is non-singular, it can be well-approximated by ordinary shape functions. The singular terms in the second term (with $j\ge 1$) are included in the local enrichment of XFEM to capture the local singular behavior at the edges.
  
\section{Numerical method}\label{sect:Numerical Method}
  Since the XFEM is an extension of the conventional FEM by including singular basis function in the shape function, we will in this section start with very brief introduction of the conventional FEM, followed by more details of XFEM as well as different local enrichment strategies for XFEM. General description of conventional FEMs can be found in many textbooks \citep[see, e.g.][]{zienkiewicz2005finite,hughes2012finite,reddy2019introduction}.
  
  \subsection{Finite Element Method}
  \par In a FEM formulation for a potential-flow problem, the fluid domain is discretized into elements, also called finite elements, and the velocity potential in each element can be approximated as
  \begin{equation}
    \label{Eq.11}
    \varphi =\sum\limits_{j=1}^{n_p}{{{N}_{j}}}(x,y){{\varphi }_{j}}.
  \end{equation}
  Here $N_j(x,y)$ is shape function, $n_p$ is the number of the nodes in the whole fluid domain and $\varphi_{j}$ denotes the nodal value of the velocity potential at node $j$. 
  Application of the Galerkin method leads to 
  \begin{equation}
  \label{Eq.12}
  \iint_{\Omega }{{{N}_{i}}}(x,y)\left[ {{\nabla }^{2}}\sum\limits_{j=1}^{n_p}{{{N}_{j}}}(x,y){{\varphi }_{j}} \right]\diff\Omega =0.
  \end{equation}
  Considering a general BVP with Dirichlet boundary $S_D$, Neumann boundary $S_N$ and Robin boundary $S_R$, the weak form of the integral in Eq.~\eqref{Eq.12} can be obtained by applying the Green's theorem and letting the test functions equal to zero on $S_D$
  \begin{equation}
  \label{Eq.13}
  	\begin{aligned}
    \int_{S_N + S_R} N_i\frac{\partial \varphi}{\partial n} \diff S
  	-\iint_{\Omega } \nabla N_i \cdot \sum\limits_{j\in S_D}\varphi_j\nabla N_j \diff\Omega \\
  	-\iint_{\Omega }{\nabla }{{N}_{i}}\cdot \sum\limits_{j\notin {{S}_{D}}}{{{\varphi }_{j}}}\nabla {{N}_{j}}\diff\Omega 
  	=0 \quad (i\notin {{S}_{D}}).
  	\end{aligned}
  \end{equation}
  For a mixed Dirichlet-Neumann BVP, as we will study in Sect.~\ref{sect:flat-plate} for the flat plate in infinite fluid, $\varphi=f_p$ on $S_D$ and $\partial \varphi  / \partial n = f_n $ on $S_N$ are known from the the boundary conditions, respectively. 
  In this case, Eq.~\eqref{Eq.13} will be represented by a linear system as
  \begin{equation}
  \label{Eq.14}
    \mathbf{K\Phi = B},
  \end{equation}
  where 
  \begin{equation}
  \label{Eq.15}
    \mathbf{\Phi} ={{\left[ \begin{matrix}
  		{{\varphi }_{1}} & {{\varphi }_{2}} & \cdots  & \begin{matrix}
  		{{\varphi }_{i}} & \cdots   \\
  		\end{matrix}  \\
  		\end{matrix} \right]}^{T}},
  \end{equation}
  Here the  superscript $T$ represents the transpose of a matrix or a vector. The elements in matrix $\mathbf{K}$ and vector $\mathbf{B}$ are defined respectively as
  \begin{equation}
  \label{Eq.16}
  %\begin{matrix}
    {{K}_{ij}} = \iint_{\Omega }{\nabla }{{N}_{i}}\cdot \nabla {{N}_{j}}\diff\Omega,  \quad (i\notin {{S}_{D}},j\notin {{S}_{D}})  %\\
  %\end{matrix}
  \end{equation}
  
  \begin{align}
  \label{Eq.17}
    {{B}_{i}}=\int_{{{S}_{b}}}{{{N}_{i}}{{f}_{n}}}\diff S
    -\iint_{\Omega }{\nabla }{{N}_{i}}\cdot \sum\limits_{j\in {{S}_{p}}}{{{({{f}_{p}})}_{j}}}\nabla {{N}_{j}}\diff\Omega, \quad (i\notin {{S}_{D}}). 
  \end{align}
  
  \par For a mixed Neumann-Robin BVP, as we will study in Sect.~\ref{sect:heave-rectangle} for the linear frequency-domain solution of a heaving rectangle at the free surface,  the weak form can be more specifically written as:
  \begin{equation}
  \label{Eq.40}
  	\begin{aligned}
  	 \iint_{\Omega }{\nabla {{N}_{i}}\cdot }\sum\limits_{j}{{{\varphi }_{j}}\nabla {{N}_{j}}}\diff\Omega 
  	 +\mathrm{i} k\int_{{{S}_{m}}}{{{N}_{i}}\sum\limits_{j}{{{\varphi }_{j}}{{N}_{j}}}}\diff S\\
  	 -k\tanh kh\int_{{{S}_{f}}}{{{N}_{i}}\sum\limits_{j}{{{\varphi }_{j}}{{N}_{j}}}}\diff S=\int_{{{S}_{b}}}{{{N}_{i}}{f_n}}\diff S.
  	\end{aligned}
  \end{equation}
  Here the mean free surface $S_f$ and control surface $S_m$ are Robin boundaries, where the boundary conditions are defined in Eqs.~\eqref{Eq.37} and \eqref{Eq.38}, respectively. The Neumann boundary condition on $S_b$ has been defined in Eq.~\eqref{Eq.36}. 
  
  \subsection{Extended Finite Element Method (XFEM)}
  XFEM was developed based on the concept of partition of unity (PU), and the so-called PU means a set of non-zero function $N_i(x,y)$  in the partition of unity domain satisfying the following condition:
  \begin{equation}
  \label{Eq.18}
    \sum\limits_{i}{{{N}_{i}(x,y)}} = 1.
  \end{equation}
  For any function in the PU domain, the following relationship holds: 
  \begin{equation}
  \label{Eq.19}
    \sum\limits_{i}{{{N}_{i}}(x,y)}\psi (x,y)=\psi (x,y).
  \end{equation}
  Undoubtedly, Eq.~\eqref{Eq.19} is also satisfied when $\psi (x,y)$ is a constant.
  Obviously, standard shape functions, for instance those shown in Eqs.~\eqref{Eq.54} and \eqref{Eq.55}, are PU functions. 
  
  After introducing the conventional FEM and the conception of PU, the enrichment function and extra degrees of freedom (DOFs) at the selected nodes will be presented. For simplicity and without losing generality, we denote $\mathcal{I}$ as the set of all nodes in the fluid domain and $\mathcal{J}$ as the subset of nodes which will be enriched. Thus the trial solution in the fluid domain with only one enrichment function on each point $j\in \mathcal{J}$ can be written as 
  \begin{equation}
  \label{Eq.20}
    \varphi =\sum\limits_{j\in \mathcal{I}}{{{N}_{j}}}(x,y){{\varphi }_{j}}+\sum\limits_{j\in \mathcal{J}}{{{N}_{j}}}(x,y)\psi (x,y){{\Psi }_{j}},
  \end{equation}
  where $\Psi_j$ represent the additional DOF at the enriched node $j$. $N_j(x,y)$ is the standard finite-element shape function, $\psi(x,y)$ denotes the enrichment function representing special knowledge, e.g. logarithmic singularity, of the fluid solution. The products ${{{N}_{j}}}(x,y)\psi (x,y)$ may be considered as local enrichment function, as their supports coincide with those of conventional finite-element shape functions, leading to sparsity in the discrete equation \citep{2010Fries}. It can be understood from Eq.~\eqref{Eq.20} that the values on nodes $j\in \mathcal{J}$ differ from $\varphi_j$, which is an unfavorable property. To ensure that nodal values are always $\varphi_{j}$ at the enriched nodes $j\in \mathcal{J}$, the enrichment function can be shifted and Eq.~\eqref{Eq.20} can be rewritten as \citep[see, e.g.][]{2010Fries, 2000Arbitrary}, 
  \begin{equation}
  \label{Eq.21}
  \begin{aligned}
    \varphi =\sum\limits_{j\in \mathcal{I}}{{{N}_{j}}}&(x,y){{\varphi }_{j}}+\sum\limits_{j\in \mathcal{J}}{{{N}_{j}}}(x,y)[\psi (x,y)-\psi (x_j,y_j)]{{\Psi }_{j}}.
  \end{aligned}
  \end{equation}
  As a result of the shifting, the enrichment represented by the 2nd summation on the right-hand side of Eq.~\eqref{Eq.21} vanishes at the nodes $j\in \mathcal{J}$, and thus recover the Kronecker-$\delta$ property of standard finite-element approximations. Unless otherwise redefined, all the enrichment functions that will be used in this paper are the shifted enrichment functions.
  
  More generally, if more than one enrichment function are introduced at each node $j\in \mathcal{J}$, Eq.\eqref{Eq.21} can be extended as  
  \begin{equation}
  \label{Eq.22}
  \begin{aligned}
    \varphi =\sum\limits_{j\in I}{{{N}_{j}}}&(x,y){{\varphi }_{j}}+\\
    &\sum\limits_{j\in J}\sum\limits_{l}{{{N}_{j}}}(x,y)[\psi^l (x,y)-\psi^l (x_j,y_j)]{{\Psi }^{l}_{j}}.
  \end{aligned}
  \end{equation}
  Here $\psi^l(x,y)$ is the $l$-th enrichment function, $\psi^l(x_j, y_j)$ denotes the value of $\psi^l(x,y)$ at $j$-th node, $\psi^l(x, y)-\psi^l(x_j, y_j)$ denotes the shifted enrichment function with a shifted value of $\psi^l(x_j, y_j)$. For brevity,  we use a matrix form to express Eq.~\eqref{Eq.22} and rewrite it as
  \begin{equation}
  \label{Eq.23}
    \varphi =\left[ \begin{matrix}
     \mathbf{N}_{std} &  \mathbf{N}_{enr}  \\
    \end{matrix} \right]\left[ \begin{matrix}
    \mathbf{\Phi}   \\
    \mathbf{\Psi}   \\
    \end{matrix} \right],
  \end{equation}
  where
  \begin{equation*}
  \begin{aligned}
    &{{N}_{stdj}}={{N}_{j}}(x,y) \quad (j=1,\cdots,n_p),\\
    &{{N}^{l}_{enrj}}={{N}_{j}}(x,y)\cdot [\psi^l (x,y)-\psi^l (x_j,y_j)] \quad \\
    &\qquad\qquad(j\in\mathcal{J}, l=1,\cdots, n^{enr}),
  \end{aligned}
  \end{equation*}
  are the elements in $\mathbf{N}_{std}$ and $\mathbf{N}_{enr}$. $n^{enr}$ denotes the number of enrichment functions. The dimension of $\mathbf{N}_{std}$ is $1\times n_p$. If there are $n^{enr}_{p}$ nodes enriched in whole domain, the dimension of $\mathbf{N}_{enr}$ is $1\times (n^{enr}_{p}\cdot n^{enr})$. Substituting Eq.~\eqref{Eq.23} into Eq.~\eqref{Eq.13}, we obtain the following expression:
  \begin{equation}
  \label{Eq.24}
  	\begin{aligned}
  	&\int_{{{S}_{N}+{S}_{R}}}{\left[ \begin{matrix}
  		\mathbf{N}_{std}^{T}  \\
  		\mathbf{N}_{enr}^{T}  \\
  		\end{matrix} \right]{\frac{\partial \varphi}{\partial n}}}\diff S-\iint_{\Omega } \nabla \mathbf{N}_{std} \cdot \mathbf{N}_{std}^D \diff\Omega \mathbf{\Phi}_D\\
  	&-\iint_{\Omega }{\left[ \begin{matrix}
  		\nabla \mathbf{N}_{std}^{T}  \\
  		\nabla \mathbf{N}_{enr}^{T}  \\
  		\end{matrix} \right]\left[ \begin{matrix}
  		\nabla \mathbf{N}_{std} & \nabla \mathbf{N}_{enr}  \\
  		\end{matrix} \right]}\diff\Omega \left[ \begin{matrix}
  	\mathbf{\Phi}   \\
  	\mathbf{\Psi}   \\
  	\end{matrix} \right]=0.
  	\end{aligned}
  \end{equation}
 Here $\mathbf{N}_{std}^D$ denotes the shape function which lies on Dirichlet boundary, and $\mathbf{\Phi}_D$ represents the velocity potential of the nodes which are located on the Dirichlet boundary. We must emphasize that there are not any enrichment nodes on the Dirichlet boundary. For a mixed Dirichlet-Neumann BVP, in the same manner as Eq.~\eqref{Eq.14}, the linear system comes from Eq.~\eqref{Eq.24} can be written as: 
  \begin{equation}
  \label{Eq.25}
    \mathbf{K X}=\mathbf{B}.
  \end{equation}
  The coefficient matrix of $\mathbf{K}$ can be divided into four parts as follows:
  \begin{equation}
  \label{Eq.26}
    \mathbf{K}=\left[ \begin{matrix}
    {\mathbf{K}^{\varphi \varphi }} & {\mathbf{K}^{\varphi \psi }}  \\
    {\mathbf{K}^{\psi \varphi }} & {\mathbf{K}^{\psi \psi }}  \\
    \end{matrix} \right]
  \end{equation}
  where the elements in $\mathbf{K}^{\varphi \varphi }$, $\mathbf{K}^{\varphi \psi }$, $\mathbf{K}^{\psi \varphi }$ and $\mathbf{K}^{\psi \psi }$ are
  \begin{equation}
  \label{Eq.27}
  \begin{aligned}
	{{K}_{ij}^{\varphi \varphi }} 
	=& \iint_{\Omega }{\nabla }{{N}_{i}}\cdot \nabla {{N}_{j}}\diff\Omega  \quad\\
	& (i\notin {{S}_{D}},j\notin {{S}_{D}}, i\in\mathcal{I}, j\in\mathcal{I}),
  \end{aligned}
  \end{equation}
  \begin{equation}
  \label{Eq.28}
  \begin{aligned}
	  {{K}_{ij}^{\varphi \psi l}} 
	=& \iint_{\Omega }{\nabla }{{N}_{i}}\cdot \nabla \left[ {{N}_{j}\left(\psi^l (x,y)-\psi^l (x_j,y_j)\right)} \right]\diff\Omega  \\
	& (i\notin {{S}_{D}},j\notin {{S}_{D}}, i\in\mathcal{I}, j\in\mathcal{J}),
  \end{aligned}
  \end{equation}
  \begin{equation}
  \label{Eq.29}
  \begin{aligned}
	{{K}_{ij}^{\psi \varphi l}} 
	=& \iint_{\Omega } \nabla \left[ {{N}_{i}\left(\psi^l (x,y)-\psi^l (x_i,y_i)\right)} \cdot {\nabla }{{N}_{j}} \right]\diff\Omega  \\
	& (i\notin {{S}_{D}},j\notin {{S}_{D}}, i\in\mathcal{J}, j\in\mathcal{I}),
  \end{aligned}
  \end{equation}
  \begin{equation}
  \label{Eq.30}
  \begin{aligned} 
	{{K}_{ij}^{\psi \psi l}} 
	=& \iint_{\Omega }\nabla \left[ {{N}_{i}\left(\psi^l (x,y)-\psi^l (x_i,y_i)\right)}\right] \cdot \\
	& \nabla \left[ {{N}_{j}\left(\psi^l (x,y)-\psi^l (x_j,y_j)\right)} \right]\diff\Omega  \\
	& (i\notin {{S}_{D}},j\notin {{S}_{D}}, i\in\mathcal{J}, j\in\mathcal{J}).
  \end{aligned}
  \end{equation}
  $\mathbf{K}^{\varphi\varphi}$ comes from conventional standard finite elements, which is only relevant to standard shape function. $\mathbf{K}^{\varphi\psi}$, $\mathbf{K}^{\psi\varphi}$ and $\mathbf{K}^{\psi\psi}$ are related to the enrichment. $\mathbf{X}$ is a $(n_p+n^{enr}_{p}\cdot n^{enr})\times1$ vector.
  The right-hand-side vector in Eq.~\eqref{Eq.25} is 
  \begin{equation}
  \label{Eq.31}
    \mathbf{B}=\left[ \begin{matrix}
    {\mathbf{B}^{\varphi}}  \\
    {\mathbf{B}^{\psi}}  \\
  \end{matrix} \right],
  \end{equation}
  where the elements in $\mathbf{B}^{\varphi}$ and $\mathbf{B}^{\psi}$ are
  \begin{equation}
  \label{Eq.32}
  \begin{aligned}
    {{B}^{\varphi }_i}
    =& \int_{{{S}_{N}+{S}_{R}}}{N_i{\frac{\partial \varphi}{\partial n}}}\diff S-\iint_{\Omega } \nabla N_i \cdot \sum\limits_{j\in S_D}\varphi_j\nabla N_j \diff\Omega\\
    & (i\notin {{S}_{D}}, i\in\mathcal{I}),
  \end{aligned}
  \end{equation}
  \begin{equation}
  \label{Eq.33}
  \begin{aligned}
    {{B}^{\psi l}_i}
    =& \int_{{{S}_{N}+{S}_{R}}}{ {{N}_{i}\left(\psi^l (x,y)-\psi^l (x_i,y_i)\right)} {\frac{\partial \varphi}{\partial n}}}\diff S.\\
    & (i\notin {{S}_{D}}, i\in\mathcal{J}).
  \end{aligned}
  \end{equation}
  $\mathbf{B}^\varphi$ is related to the standard FEM and $\mathbf{B}^\psi$ is related to the local enrichment. To be clear, Eq.~\eqref{Eq.27} is equivalent to Eq.~\eqref{Eq.16}, and Eq.~\eqref{Eq.32} is equivalent to Eq.~\eqref{Eq.17} in the previous section.
  
  \par For a mixed Neumann-Robin BVP, we use corner-flow solution as the enrichment solution and Eq.~\eqref{Eq.22} to construct the local approximation, and obtain a final equation system similar to Eq.~\eqref{Eq.40} as:
  \begin{equation}
  	\label{Eq.41}
  	\begin{aligned}
  	&\iint_{\Omega }{\left[ \begin{matrix}
  		\nabla \mathbf{N}_{std}^{T}  \\
  		\nabla \mathbf{N}_{enr}^{T}  \\
  		\end{matrix} \right]\left[ \begin{matrix}
  		\nabla {\mathbf{N}_{std}} \nabla {\mathbf{N}_{enr}}  \\
  		\end{matrix} \right]}\diff\Omega \left[ \begin{matrix}
  	\mathbf{\Phi}   \\
  	\mathbf{\Psi}   \\
  	\end{matrix} \right]\\
  	&+\ii k\int_{{{S}_{m}}}{\mathbf{N}_{std}^{T}} \cdot {\mathbf{N}_{std}}\diff S\cdot \mathbf{\Phi} -\\
  	& k\tanh kh \int_{{{S}_{f}}}{\mathbf{N}_{std}^{T}\cdot{\mathbf{N}_{std}}}\diff S\cdot \mathbf{\Phi}
  	=\int_{{{S}_{b}}}{\left[ \begin{matrix}
  		\mathbf{N}_{std}^{T}  \\
  		\mathbf{N}_{enr}^{T}  \\
  		\end{matrix} \right]{{f}_{n}}}\diff S. 
  	\end{aligned}
  \end{equation}
  \begin{figure}[t]
    	\centering
    	\includegraphics[scale=.4]{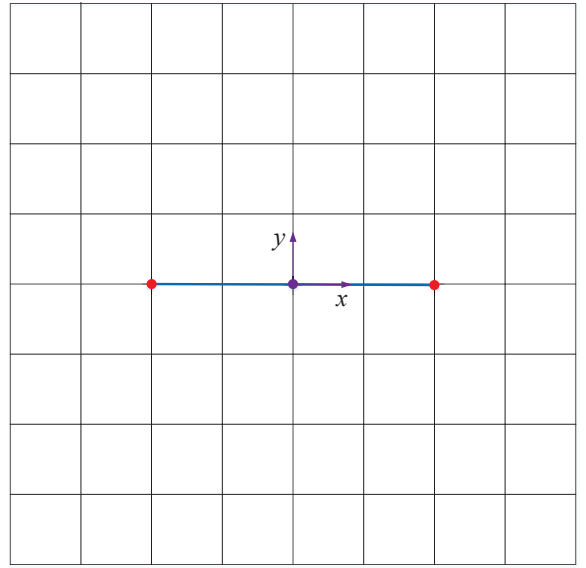}
    	\caption{An illustration of the point enrichment.}
    	\label{Fig.2}	
    \end{figure}
  \begin{figure}[t]
   	\centering
   	\includegraphics[scale=.4]{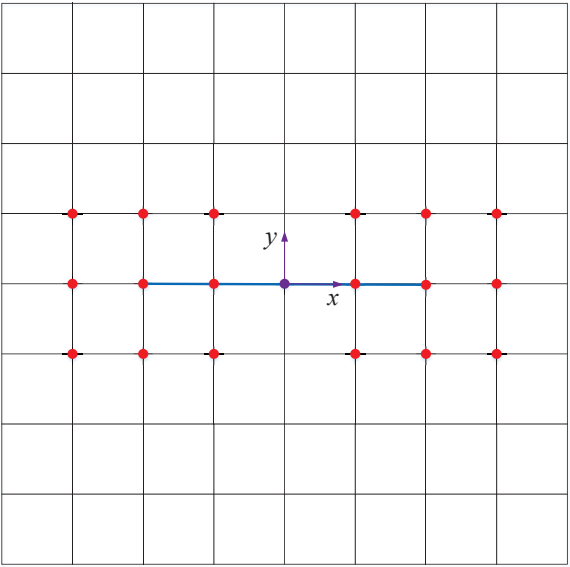}
   	\caption{An illustration of the patch enrichment.}
   	\label{Fig.3}	
   \end{figure}
   \begin{figure}[t]
   	\centering
   	\includegraphics[scale=.4]{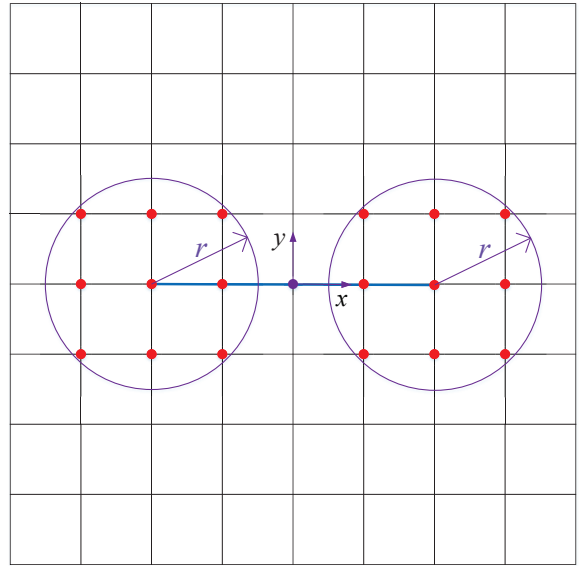}
   	\caption{An illustration of the radius enrichment.}
   	\label{Fig.4}	
   \end{figure}
  \subsection{Enrichment strategies}\label{sect:enrichments}
   \par In the previous subsection, the XFEM has been introduced through mixed BVPs. The key point of XFEM is the local enrichment, and we will in this subsection discuss three different enrichment strategies in detail. Unless otherwise mentioned in present work, a \textit{singular point} is defined as a point where the fluid velocity is infinite, and an element is called \textit{singular element} if it contains at least one \textit{singular point}. A \textit{singular patch} is a patch of multiple elements, among which at least one is a \textit{singular element}.
   
   \subsubsection{Point enrichment}\label{sect:point-enri}
    \par In the point enrichment approach, as depicted in Fig.~\ref{Fig.2}, singular solutions are introduced to enrich the local approximation only on the \textit{singular points}, and the end points of the blue line are the \textit{singular points}. In this way, the additional number of unknowns due to enrichment is only dependent on the number of \textit{singular points} and the number of singular terms introduced at each \textit{singular point}, and thus is independent on the meshes. Therefore, this enrichment  only influences the \textit{singular elements} which contain the \textit{singular points}. The influence domain of a enriched point depends on the mesh size. Consequently, the enhancement of the solution accuracy may not increase as the mesh are refined, which will be discussed later in Sect.~\ref{sect:flat-plate}.
    
   \subsubsection{Patch enrichment}\label{sect:patch-enri}
   \par Compared with point enrichment, the patch enrichment method will introduce enrichment to all points on the \textit{singular patch}, which are represented by the filled circles in Fig.~\ref{Fig.3}. The patch enrichment includes the first layer of neighboring points surrounding the singular points. Similar to the point enrichment method, the enrichment domain of patch enrichment depends on the mesh size, and the additional number of unknowns do not increase even if the meshes are refined. Similar to the point enrichment, patch enrichment experiences low convergent rate with the refinement of the local meshes.
   
   \subsubsection{Radius enrichment}\label{sect:radius-enri}
   \par Different from the point enrichment and patch enrichment, the radius enrichment method will enriches the solution at all points within a circle with predefined radius $R_{enri}$. Here $R_{enri}$ must be a positive value and it is independent of the mesh size. As demonstrated in Fig.~\ref{Fig.4}, the center point of the enrichment area is the \textit{singular point}. The value of $R_{enri}$ may be taken as $1/10$ of the space dimension, as it was suggested by \citep{laborde2005high}. In the present work, we normally take $R_{enri}=0.2$ as we are considering 2D problems. More detail about how to choose the enrichment radius will be discussed in the numerical example of heaving rectangle at free surface. The drawback of this enrichment method is that the additional number of unknowns will increase with the mesh refinements.
   
   \subsection{Integrals on the singular elements}\label{sect:singular integration}
    In this subsection, the integration on singular elements is discussed. For illustrating the element integral strategy explicitly, a singular function stems from corner-flow solution in Sect.~\ref{sect:corner-flow} is utilized as the enrichment function. In Eq.~\eqref{Eq.9}, the most singular term is the first term, i.e. the term with $j=1$ 
    \begin{equation}
    	\label{Eq.60}
    	\psi(x,y) = r^{m_1}\cos(m_1\theta).
    \end{equation}
    For demonstration purpose, we will take this term as an example of the enrichment function, and discuss numerical integration of singular terms on the elements. In practice, more terms in Eq.~\eqref{Eq.9} can be included as enrichment functions, using the similar procedure that will be described in the rest of this section. 
    
    The trial solutions in Eqs.~\eqref{Eq.20} and \eqref{Eq.21} involve the evaluation of the following enrichment shape function 
    \begin{equation}
    	\label{Eq.61}
    	N^{enr}_j= N_{j}(x,y) \psi(x,y) = N_{j}(x,y) r^{m_1}\cos(m_1\theta).
    \end{equation}
    Here $(x,y)$ is the location in physical space, which can be obtained from isoparametric element illustrated in Fig.~\ref{Fig.22}. Derivatives of the enrichment shape  function with respect to $x$ and $y$ are expressed by 
    \begin{equation}
    \label{Eq.62}
    	\begin{aligned}
    	& \frac{\partial {N^{enr}_j}}{\partial x}=\frac{\partial {{N}_{j}}}{\partial x}{{r}^{m_1}}\cos (m_1\theta )+{{N}_{j}}\frac{\partial }{\partial x}\left( {{r}^{m_1}}\cos (m_1\theta ) \right), \\ 
    	& \frac{\partial {N^{enr}_j}}{\partial y}=\frac{\partial {{N}_{j}}}{\partial y}{{r}^{m_1}}\cos (m_1\theta )+{{N}_{j}}\frac{\partial }{\partial y}\left( {{r}^{m_1}}\cos (m_1\theta ) \right). \\ 
    	\end{aligned}
    \end{equation}
    Substituting Eq.~\eqref{Eq.62} into Eq.~\eqref{Eq.30}, the diagonal entry of enriched element stiffness matrix ${K}_{ii}^{\psi \psi}$ can be written as
    \begin{equation}
    	\label{Eq.63}
    	K_{ii}^{\psi \psi }=\int_{{{\Omega }^{e}}}{{{\left( \frac{\partial {N^{enr}_i}}{\partial x} \right)}^{2}}+}{{\left( \frac{\partial {N^{enr}_i}}{\partial y} \right)}^{2}}\diff s,
    \end{equation}
    where $\Omega^e$ denotes the surface of elements.  Point $i$ is one of their nodes. Apparently, if the interior angle $\gamma < \pi$, the $x$- and $y$-derivatives of $ r^{m_1}\cos (m_1\theta )$ are singular, with a singularity of $r^{m_1-1}$ as $r\rightarrow 0$. Thus the square terms in Eq.~\eqref{Eq.63} are $r^{2m_1-2}$ singularities. It is challenging but important to accurately calculate this singular integration. In this paper, the so-called DECUHR adaptive quadrature algorithm \citep{1994DECUHR} is employed to overcome the difficulties in numerical integration of Eq.~\eqref{Eq.63}. The DECUHR algorithm combines an adaptive subdivision strategy with an extrapolation of the error expansion, where a non-uniform subdivision of the element close to a \textit{singular point} is employed. More details of the DECUHR algorithm can be found in \cite{1994DECUHR}, and the application of this algorithm in GFEM to deal with singular integrals can be found in \cite{2000Thedesign}. An open-source FORTRAN code of this DECUHR algorithm from Alan Genz Software website of Washington State University, which can deal with problems at the dimension of 2$\sim$15, has been applied in this study. It is not an option in the code to handle 1D singular integrals. 
    
    \par In present work, an adaptive Gaussian quadrature algorithm is applied to accurately calculate the 1D singular integrals. It consists of the following steps:
\begin{steps}
\item Setting a fixed \textit{tolerance}, using $T$ to represent the result, letting $T=0$, using $T_3$ represent the temporary variate and $T_3=0$. 
\item Using Gaussian integral to obtain the numerical integration result in whole element and written as $T_1$, letting $T=T_1$.
\item Dividing the element into two uniform element, and integrating in those two subdivision element, expressing as $T_{21}$ and  $T_{22}$ respectively, assuming $T_{21}$ is the numerical result in singular element, letting $T_{2} =T_{21}+T_{22}$ and $T_3 =T_3+T_{22}$ . 
\item Calculate $error = |T_2 -T_1|$. If $error > tolerance$, we denote $T =T_3+T_{21}$, divide the sub-element which contains \textit{singular point} into two, and go to Step 2. Otherwise, if $error \le tolerance$, we output $T$ as the final result.
\end{steps}

   \begin{figure}[t]
  	\centering
  	\includegraphics[scale=.413]{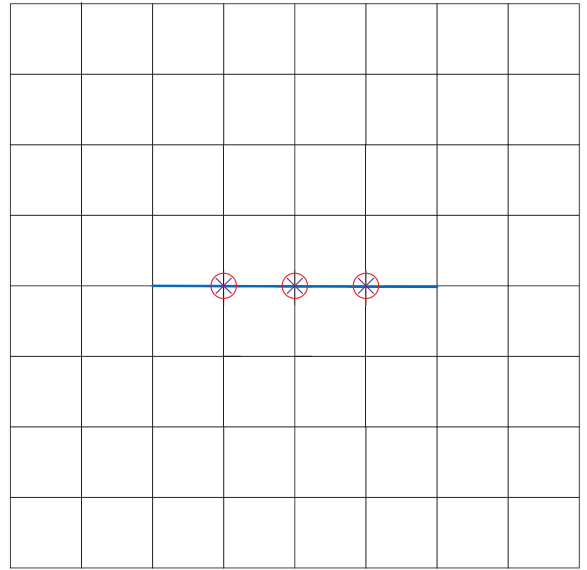}
  	\caption{An illustration of the double nodes on a flat plate. The blue line represents the flat plate with vanishing thickness.}
  	\label{Fig.5}	
  \end{figure}
  \begin{figure*}[t]
  	\centering
  	\subfigure[Velocity potential in the fluid] {\includegraphics[scale=.30]{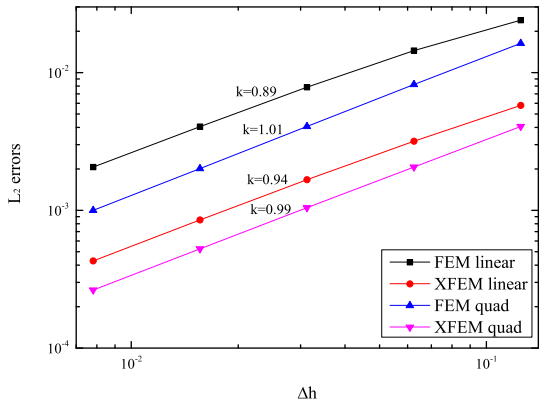}}
  	\subfigure[Added mass] {\includegraphics[scale=.30]{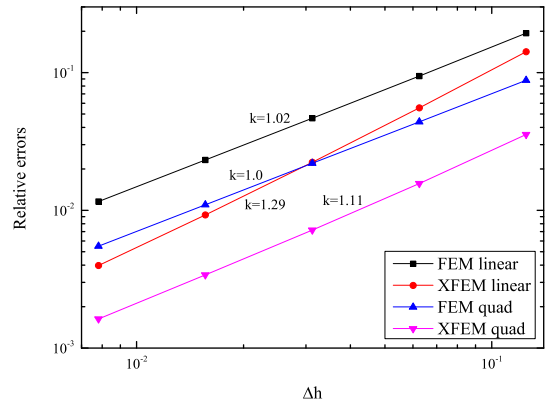}}
  	\caption{Results of mesh-convergence study for four FEMs using the point enrichment approach. $\Delta h$ = mesh size, $k$=slope.}
  	\label{Fig.6}	
  \end{figure*}
  
\section{Numerical studies}\label{sect:Numerical case}
  \par For verification purposes, an uniform flow around a flat plate of vanishing thickness in 2D, and the added mass of the same plate are considered. Then, the heaving rectangle on a free surface is studied via linear FEM, linear XFEM, quadratic FEM and quadratic XFEM.
  
  \subsection{Uniform flow around a flat plate and added mass of a flat plate}\label{sect:flat-plate}
  \begin{figure*}[t]
  	\centering
  	\subfigure[Velocity potential in the fluid] {\includegraphics[scale=.30]{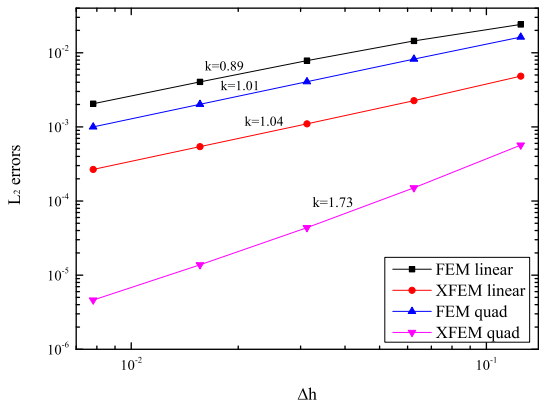}}
  	\subfigure[Added mass] {\includegraphics[scale=.30]{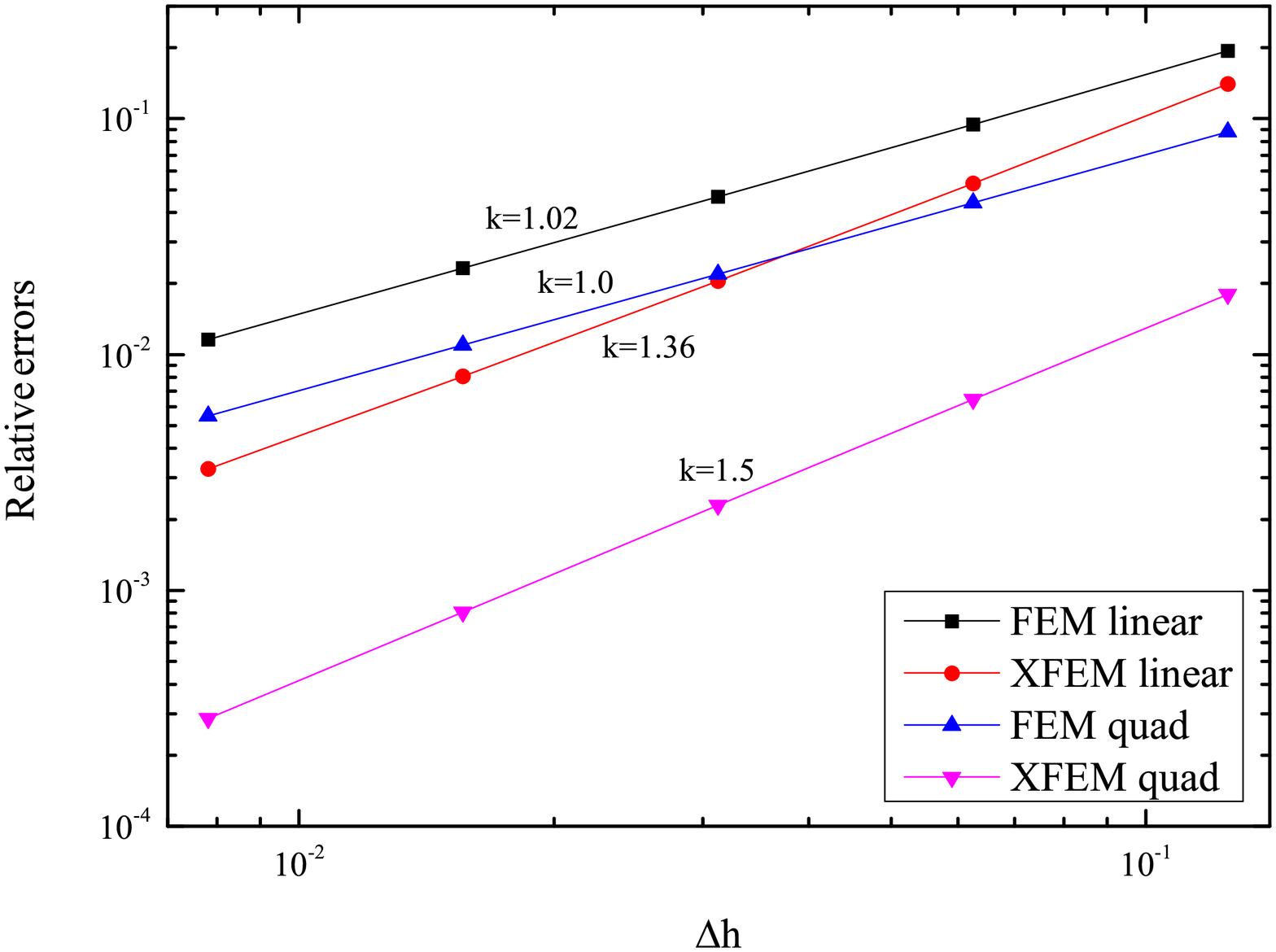}}
  	\caption{Results of mesh-convergence study for four FEMs using the patch enrichment approach. $\Delta h$ = mesh size, $k$=slope.}
  	\label{Fig.8}	
  \end{figure*}
  
  \begin{figure*}[t]
  	\centering
  	\subfigure[Velocity potential in the fluid] {\includegraphics[scale=.30]{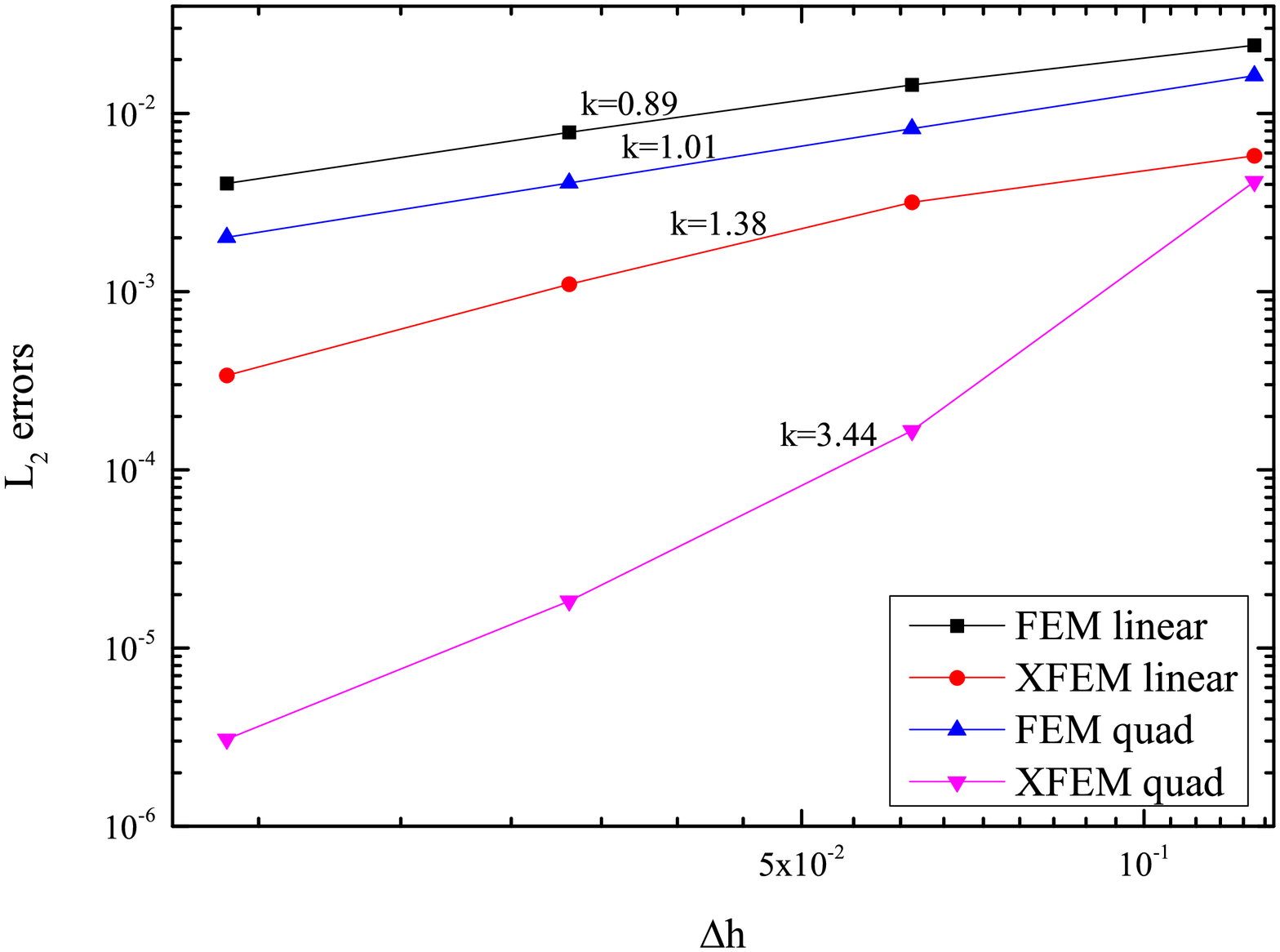}}
  	\subfigure[Added mass] {\includegraphics[scale=.30]{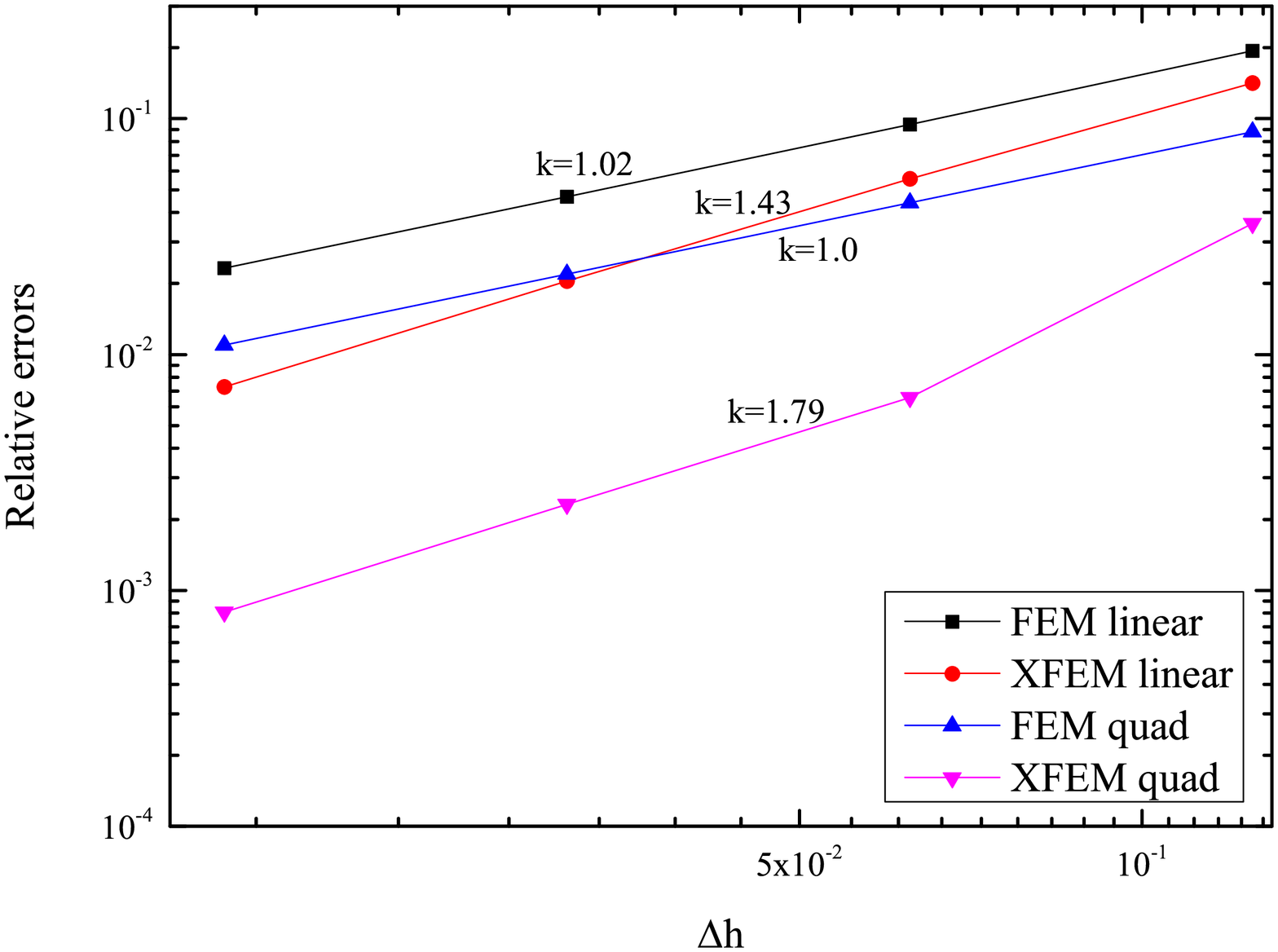}}
  	\caption{Results of mesh-convergence study for four FEMs using the radius enrichment approach. $\Delta h$ = mesh size, $k$=slope.}
  	\label{Fig.9}	
  \end{figure*}
  
  \begin{figure}[t]
  	\centering \includegraphics[scale=.30]{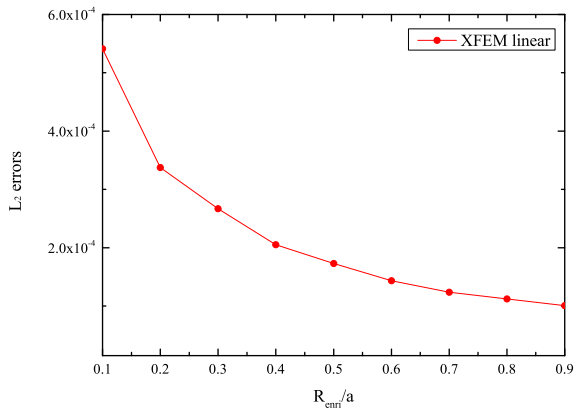}
  	\caption{The error of velocity potential versus the non-dimensional enrichment radius $R_{enri}/a$ for linear XFEM. $R_{enri}$ represents the enrichment radius, $a$ is half breadth of the plate. $64\times 64$ uniform meshes have been used.}
  	\label{Fig.10}	
  \end{figure}
  
  \begin{figure*}[t]
  	\centering
  	\subfigure[Linear methods] {\includegraphics[scale=.30]{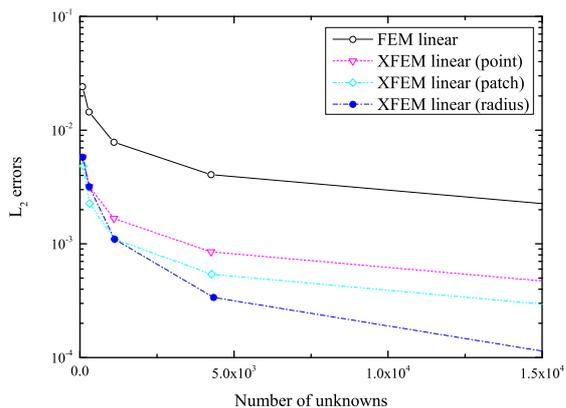}}
	\subfigure[Quadratic methods] {\includegraphics[scale=.30]{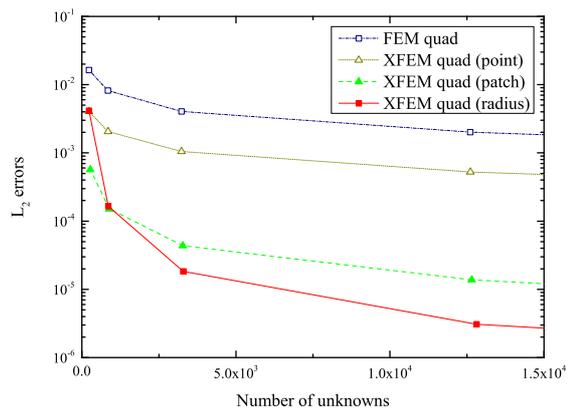}}
  	\caption{The $L_2$ errors as function of number of unknowns for the conventional FEMs and their corresponding XFEMs using different enrichment strategies.}
  	\label{Fig.11}	
  \end{figure*}
  
  \begin{figure}[t]
  	\centering
  	\includegraphics[scale=.30]{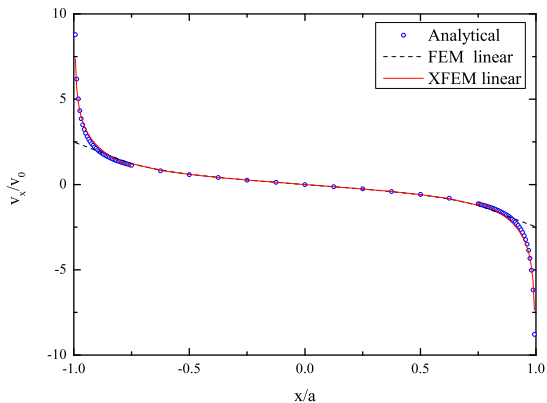}
  	\caption{Horizontal velocity distribution along the length of plate. Results are presented for conventional linear FEM and linear XFEM, together with the analytical solutions. $a$ = half breadth of the flat plate, $v_0$ = free-stream inflow velocity.}
  	\label{Fig.7}	
  \end{figure}
  \begin{table*}[htb]
  	\centering
  	\caption{The number of unknowns of linear FEM and linear XFEM with three different enrichment strategies. Four different mesh densities are considered.}
  	\begin{tabular}{ccccc}
  		\toprule
  		\multirow{2}*{Mesh size $(\Delta h/a)$} &\multirow{2}*{Linear FEM}&Linear XFEM&Linear XFEM& Linear XFEM \\
  		 &                       & (point enrichment) & (patch enrichment) & (radius enrichment) \\
  		\midrule
  		0.5	   & 84   &	86   &	104	& 86    \\
  		0.25   & 296  &	298  &	316	& 298   \\
  		0.125  & 1104 &	1106 & 1124	& 1124  \\
  		0.0625 & 4256 & 4258 & 4276 & 4336  \\
  		\bottomrule
  	\end{tabular}
  	\label{tab1}
  \end{table*}
  \begin{table*}[htb]
  	\centering
  	\caption{The number of unknowns of quadratic FEM and quadratic XFEM with three different enrichment strategies. Four different mesh densities are considered.}
  	\begin{tabular}{ccccc}
  		\toprule
	  	\multirow{2}*{Mesh size $(\Delta h/a)$} & \multirow{2}*{Quad. FEM}  & Quad. XFEM & Quad. XFEM & Quad. XFEM\\
	  	&  & (point enrichment) & (patch enrichment) & (radius enrichment) \\
  		\midrule
  		0.5    & 232   & 234   & 278   & 234   \\
  		0.25   & 848   & 850   & 894   & 860   \\
  		0.125  & 3232  & 3234  & 3278  & 3288  \\
  		0.0625 & 12608 & 12610 & 12654 & 12814 \\
  		\bottomrule
  	\end{tabular}
  	\label{tab2}
  \end{table*}

  \par The analytical solution of the complex potential for uniform flow around a 2D thin flat plate in an infinite domain can be found in the textbook of \cite{Newman2017Marine} and the Appendix B, where we also show that a modification of the sign in the original formula is needed for the flow variable on the right-half plane.  
  
  \par To model the flat plate in infinite fluid, we have to use a truncated fluid domain in our numerical method. See a sketch of the truncated domain in Fig.~\ref{Fig.5}. Based on the analytical solution, Dirichlet boundary conditions are specified at the truncated boundaries surrounding the fluid domain, and Neumann boundary condition on the upper and lower surfaces of the thin plate. The mathematical formulation of the mixed BVP has been  described in Sect.~\ref{sect:mathematical-formu} and the conventional FEM and XFEM explained in Sect.~\ref{sect:Numerical Method}. Even though the flat plate has zero thickness, the velocity potentials are different on the two sides of the plate. Thus a double-node technique is used on the plate except at the two endpoints of the flat plate, where the velocity potential must be continuous. The double-node technique allows for two velocity potential values at the same location. See an illustration in Fig.~\ref{Fig.5}, where the open circles and crosses represent two different nodes respectively.
  
  \par To solve the mixed Dirichlet-Neumann BVP numerically, we have implemented four different FEM solvers, including linear FEM, linear XFEM, quadratic FEM and quadratic XFEM. In the linear and quadratic XFEMs, we have used the analytical solution as the enrichment function at the enrichment nodes close to the singular points, in this case the two ends of the plate. Figs.~\ref{Fig.2}, \ref{Fig.3} and \ref{Fig.4} illustrate the point, patch and radius enrichment strategies, respectively. In order to compare the global accuracy of different methods, the $L_2$ errors of the velocity potential on all grid points will be presented as a function of mesh size $\Delta h = \Delta x = \Delta y$. The $L_2$ error is defined as
  \begin{equation}
  \label{Eq.34}
    e_{L_2}=\sqrt{\left.\sum_{i=1}^{N}\left(\phi_{i}^{\rm num}-\phi_{i}^{\rm ana}\right)^{2} \middle/ \sum_{i=1}^{N}\left(\phi_{i}^{\rm ana}\right)^{2}\right.},
  \end{equation}
  where $\phi_{i}^{\rm num}$ denotes numerical solution on the $i$th node, and $\phi_{i}^{\rm ana}$ represents the corresponding analytical solution. $N$ denotes the total number of nodes. 
  
  In principle, the added mass of the plate should be calculated by considering oscillatory motions of the plate. However, for the special case of a structure in unbounded fluid or practically sufficient away from any other boundaries, the added mass can also be obtained based on the flow solution for fixed structure in a uniform flow. According to Section 4.10 in \cite{Newman2017Marine}, the velocity-potential solution around a fixed structure in a free stream along $y$-axis can be expressed as $\phi_{\rm fix} = Uy - \phi_{\rm move}$, where $\phi_{\rm move} = U \phi_2$ is the velocity potential due to the same structure moving along $y$-axis with velocity $U$. $\phi_2$ is the velocity potential induced by the structure at a unit velocity, i.e. $U=1$, along $y$-axis. Therefore, $\phi_2$ can be immediately obtained as $\phi_2 = (U\,y - \phi_{\rm fix})\,/\,U$, and thus the added mass can be calculated according to Eq.(114) in \cite{Newman2017Marine}, which only needs the flow solution of $\phi_2$. Note that the above discussions are still valid if the velocity $U$ is not a constant in time, i.e. $U=U(t)$. 
  
  \par The $L_2$ errors for velocity potential and the relatively errors for added mass of different method via different enrichment strategies are presented in Fig.~\ref{Fig.6}, \ref{Fig.8} and \ref{Fig.9}, respectively.

  \par As shown in Fig.~\ref{Fig.6}(a) and \ref{Fig.6}(b), when the point enrichment as described in Sect.~\ref{sect:point-enri} is applied at the two end-points of the plate, the errors are greatly reduced, indicating that the local enrichment in XFEMs is very effective in reducing both the local and global errors, which is expected. However, it is also surprising that all FEMs, including the quadratic FEMs, showed convergent rates close to $1.0$ for the velocity potential on fluid points. The $k$ values in the figures are the fitted convergence rate based on five different mesh densities. Similarly, as seen in Fig.~\ref{Fig.6}(b), lower than expected mesh-convergence rates are observed for the added mass of the flat plate. The influence area of the enrichment functions is smaller for a locally finer mesh close to the \textit{singular point}, due to the fact that the interpolations by the finite-element shape functions in Eqs.~\eqref{Eq.21} and \eqref{Eq.22} are piecewise. The shape function $N_j$ at the point $j$ is always zero over the elements, which do not own point $j$ as one of their element nodes. At the non-enriched points, sufficiently close to the \textit{singular point} but belong to none of the \textit{singular elements}, the velocity potential also changes dramatically, and thus the applied smooth shape functions have difficulties to accurately capture the strong local singular solutions. Actually, according to \cite{zienkiewicz2005finite}, the singularity affects not only the local area, but also at a distance surrounding the singularity, and consequently the convergence rate for the ordinary FEMs. The affected convergence rate follows $O(N_{dof})^{(-min[\lambda,p]/2)}$, where $\lambda$ is a number associated with the intensity of the singularity, $p$ represents power exponent of the basis function, $N_{dof}$ the number of freedom. More details about the effect can be found in \citet[pp. 458-459]{zienkiewicz2005finite}.
  
  \par The results of convergence studies are presented in Fig.~\ref{Fig.8} for the patch enrichment. For both the solutions of velocity potential at grid nodes and added mass, only marginal increases of the convergence rate of linear XFEM are seen. However, the improvement in the results of quadratic XFEM is notable for both velocity potential and added mass. Theoretically, we expect the convergence rate of a quadratic method be equal or greater than 2. Even though the overall accuracy of linear and quadratic XFEM has been greatly improved with the adoption of patch enrichment instead of the point enrichment, it is still below our expectation, in particular for the quadratic XFEM. The reason is as follows: similar to the point enrichment strategy, the enrichment area of the patch enrichment strategy is also mesh-dependent. 
  
  \par To eliminate the mesh-dependency of the local enrichment, the radius enrichment method, as illustrated in Fig.~\ref{Fig.4}, appears to be a good choice. In this method, the enrichment area is a predefined constant and independent on mesh sizes. As demonstrated by the results of mesh-convergence study in Fig.~\ref{Fig.9}, the superiority of XFEMs, in particular the quadratic XFEM, is remarkable, when the radius enrichment is applied. We have used an constant enrichment radius of $R_{enri} = 0.2$, as suggested by \cite{laborde2005high} for 2D problems, at both ends of the plate. Compared to the conventional linear FEM, convergence rate of linear XFEM for the velocity potential advanced notably from $k=0.89$ to $k=1.38$, and from $k=1.02$ to $k=1.43$ for added mass. The convergence rate of quadratic XFEM improved exceedingly from $k=0.99$ to $k=3.44$ for velocity potential, and for added mass from $k=1.11$ to $k=1.79$.

  \par For the present case, since we are using the analytical solution as the enrichment function at the singular points, the accuracy of the XFEM solutions will further improve if a larger enrichment radius is applied. This is illustrated in Fig.~\ref{Fig.10}, where we have presented the $L_2$ errors of the velocity potential as function of $R_{enri}/a$, the ratio between enrichment radius and length of the plate. 
  
  \par For both linear and quadratic XFEMs, it is apparent from the comparisons in Figs.~\ref{Fig.6}-\ref{Fig.9} that, the radius enrichment strategy yields the most accurate results for a given mesh resolution, with a cost of introducing more extra DOFs (or unknowns in the final linear system) than the two other enrichment strategies. Since all points within a radius $R_{enri}$ to the singular points will be enriched, too many extra DOFs may be introduced if an unnecessarily large $R_{enri}$ has been chosen. For a given $R_{enri}$, the number of extra DOFs is also larger for a finer mesh. On the other hand, the point enrichment method introduces fewest extra DOFs, but the accuracy is the lowest among the three enrichment methods. From practical application point of view, it is recommended to apply the radius enrichment method with a small enrichment area.    
  
  \par It is of more interest to compare the computational efforts to achieve a similar accuracy. In this regard, we have also plotted in Fig.~\ref{Fig.11} the $L_2$ errors of the velocity potential as function of the total number of unknowns, which is an indicator of CPU time. The number of unknowns for different enrichment strategies and different mesh sizes are listed in Table~\ref{tab1} for linear FEMs and Table~\ref{tab2} for quadratic FEMs. It is apparent that the local enrichment increases only marginally the total number of unknowns, while reducing the global errors significantly.
  
  \par To illustrate how the XFEMs has increased the accuracy of the local flows, the horizontal velocity along the flat plate is shown in Fig.~\ref{Fig.7}. Here the conventional linear FEM and linear XFEM are compared. Solid line represents result for linear XFEM, while result of conventional linear FEM is represented by the dash line. The corresponding analytical solution is denoted by open circles. Thanks to the local enrichment, linear XFEM shows very encouraging results, especially close to the singular point. On the contrary, the conventional linear FEM fails to capture the strong variation of the flow at two ends of the plate. Since the applied FEM is only $C^0$ continuous, the presented velocity in the figure at each point is the average value of the velocities calculated in the elements sharing the point. For the \textit{singular elements}, the velocity was obtained by differentiating the shape functions in Eq.~\eqref{Eq.21}, and we have added more points within the element to better illustrate the variation of velocity therein. In theory, the nodal FEMs based on shape functions are only $C^0$ continuous, and thus the velocities are discontinuous between elements for both ordinary FEMs and XFEMs. In Fig.~\ref{Fig.7}, the velocity distribution appears smoother for FEM because we have used the average of the velocities evaluated on the adjacent elements as the nodal values. If the averaging is not applied, the velocities will  appear discontinuous at all nodes. Since the solution representation is more accurate in the enriched element than that in its neighboring ordinary element, a more obvious jump of the velocity has been observed (at $x/a=\pm 0.75$ in the present case).
  
  \begin{figure*}[t]
  	\centering
  	\includegraphics[scale=.55]{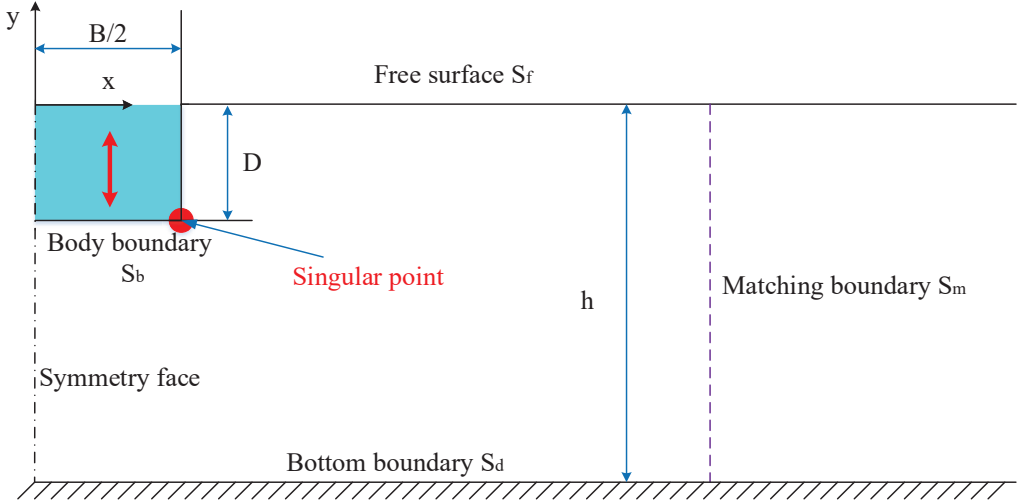}
  	\caption{Sketch map of half rectangular heaving on the free surface.}
  	\label{Fig.12}	
  \end{figure*}
 
  \subsection{Heaving rectangular cylinder on free surface}\label{sect:heave-rectangle}
  \begin{figure}[t]
  	\centering
  	\includegraphics[scale=.30]{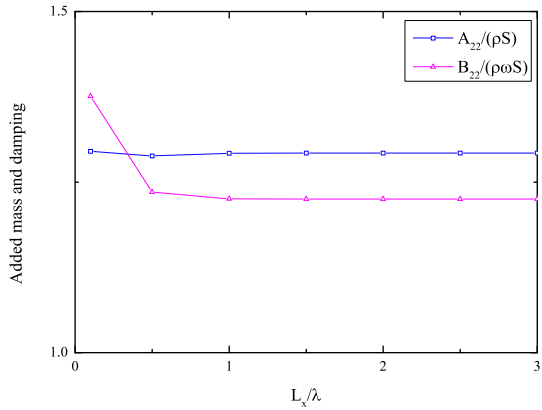}
  	\caption{Convergence performance of the horizontal length from the rectangle to the matching boundary when the square of forcing frequency is $\omega^2B/(2g)=0.1$. $A_{33}$=heave added mass, $B_{33}$=heave radiation damping, $S$=submerged cross-sectional area, $\rho$=mass density of water, $\omega$=circular frequency, $L_x$ = horizontal length from the rectangle to the matching boundary, $\lambda$ = wavelength of radiated waves. }
  	\label{Fig.13}	
  \end{figure}
  
  \begin{figure}[t]
  	\centering
  	\includegraphics[scale=.4]{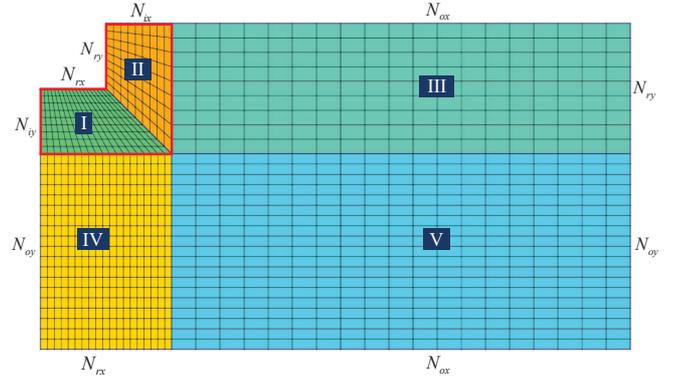}
  	\caption{Schematic of the mesh of linear elements for the half rectangular heaving on the free surface.}
  	\label{Fig.14}	
  \end{figure}
  
  \begin{figure*}[t]
  	\centering
  	\subfigure[Added mass]{\includegraphics[scale=.30]{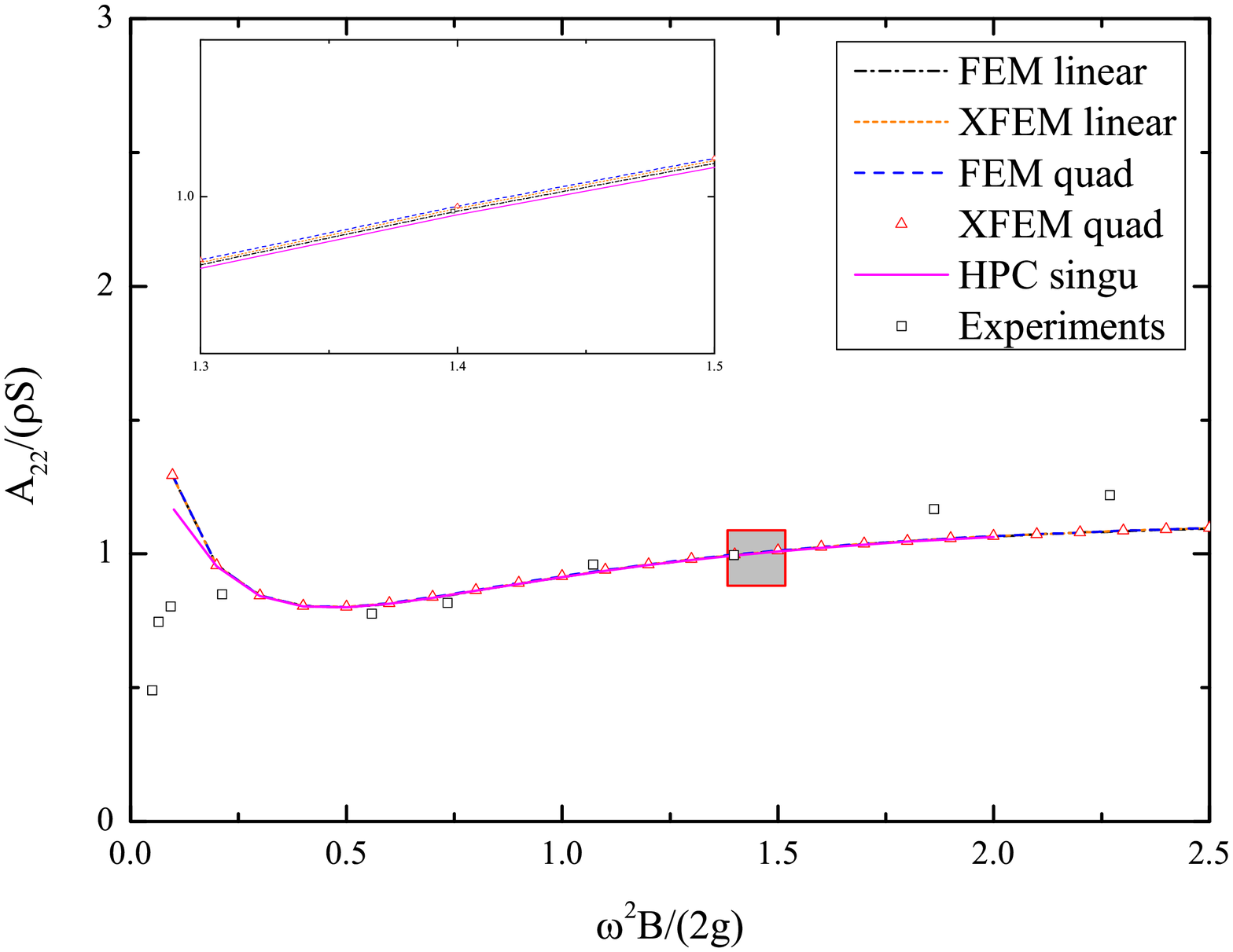}}
  	\subfigure[Radiation damping]{\includegraphics[scale=.30]{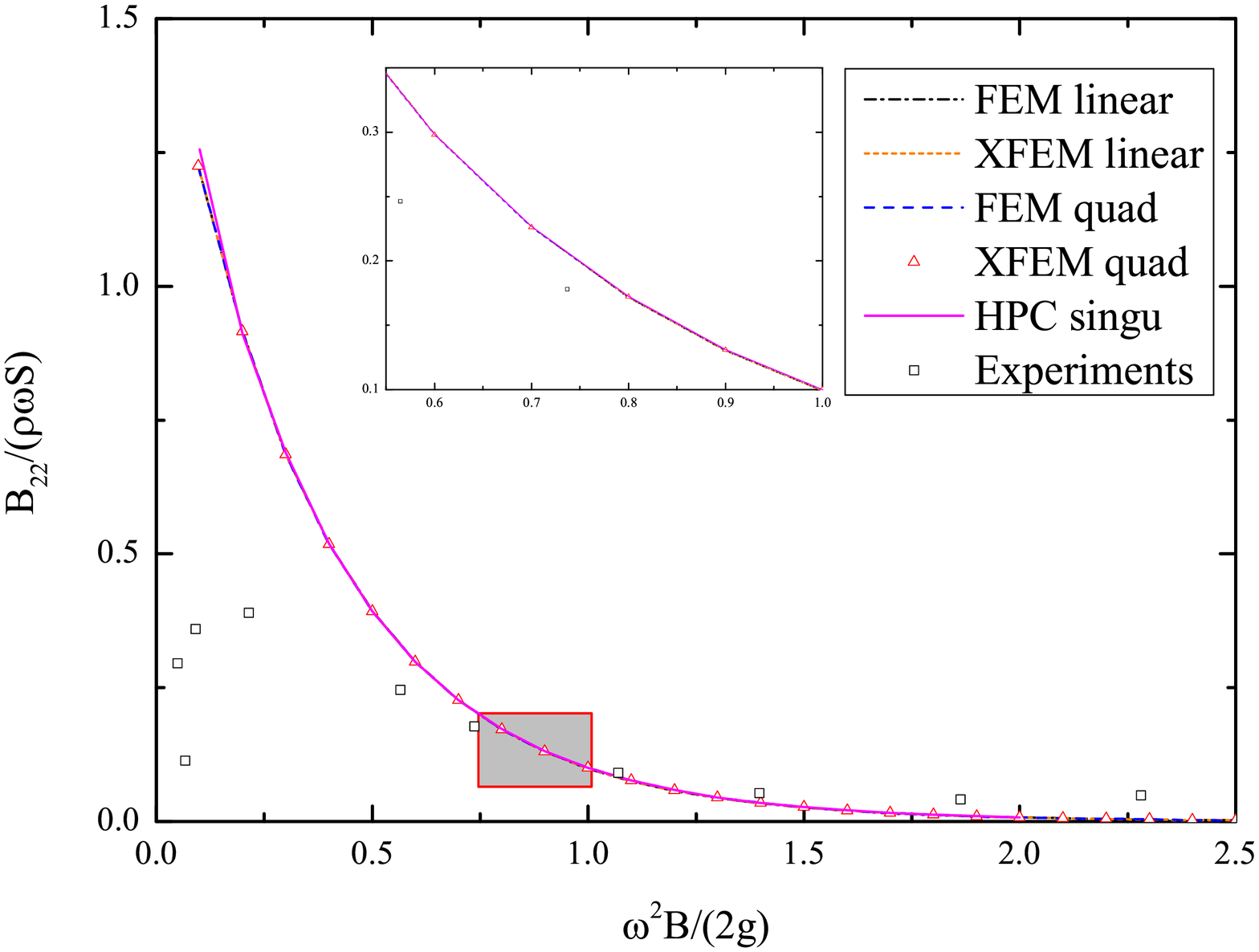}}
  	\caption{The added mass and radiation damping as function of the oscillatory frequency of a floating rectangular with beam-to-radio $(B/D)$ equals to 2.0. $B$=beam, $D$=draft, $A_{33}$=heave added mass, $B_{33}$= heave radiation damping, $S$=submerged cross-sectional area, $\rho$=mass density of water, $\omega$=circular frequency.}
  	\label{Fig.16}	
  \end{figure*}
  
  \par In this part, under the framework of linear potential-flow theory in the frequency domain, a floating heaving rectangular cylinder on a free surface is considered. $B$ and $D$ are used to represent beam and draft of the rectangle, respectively. The considered water depth is $h = 40D$, and the beam-to-draft ratio $B/D$ is taken as 2.0. An illustration of half of the domain is presented in Fig.~\ref{Fig.12}. In theory, a radiation condition should be applied at $x \rightarrow \pm \infty$. In practice, it is impossible to model an fluid domain with infinite extension, and a truncation at a certain horizontal distance $L_x$ from the rectangle must be made. The radiation condition is then applied on a control surface $S_m$, which is chosen sufficiently far from the structure. In this study, we choose a horizontal truncation distance as twice the longest wavelength that will be studied, and use the same computational domain for all different cases.  
  
  \par To reduce the computational costs, the symmetric property of the considered problem is utilized, and thus only half of the fluid domain is considered. At the symmetry plane, horizontal velocity is equal to zero, i.e. ${\partial \phi} / {\partial x} = 0$. An illustration of the computational domain is presented in Fig.~\ref{Fig.12}.
  
  \subsubsection{Linear added mass and damping coefficients}
  
  \par Fig.~\ref{Fig.13} displays the non-dimensional added mass and damping coefficients for different truncation distances $L_x$ from the rectangle. A non-dimensional wave number $\omega^2B/(2g)$ = 0.1 has been considered in the calculations, corresponding to the longest wave that will be considered in this section. If the selected $L_x$ has negligible results for the longest wave, it is also considered as sufficient for the shorter waves. It is apparent from the results in Fig.~\ref{Fig.13} that the hydrodynamic coefficients do not change as long as $L_x/\lambda \ge 1.0$. Here $\lambda$ is the wavelength. $L_x=2\lambda$ will be applied in our later analyses in this section.
  
  \par Matched multi-block meshes in the fluid domain are utilized as a starting point, with block I and block II fitted to the body surface, block IV below the body surface, block III and block V away from the structure. See an example of the meshes in Fig.~\ref{Fig.14}, generated from the open source mesh generator GMSH. The following parameters are defined to denote the number of elements along the sides of the blocks to control the mesh densities in different blocks: $N_{rx}$ is the number of elements on the bottom of the rectangle, $N_{ry}$ along the side wall of the rectangle, $N_{ix}$ along the free surface in the inner block, $N_{iy}$ in vertical direction of internal block at symmetry face. Correspondingly, $N_{oy}$ represent the element number in the horizontal direction of external domain on free surface, $N_{oy}$ represent the element number in vertical direction of external domain at symmetry face. Here $N_{rx}$ must be equal to $N_{ix}$ so that blocks I and II will match at their common boundary. For simplicity, we will also take  $N_{rx}=N_{ry}=N_{i
 x}=N_{iy}$.  Meshes in blocks IV and V  are stretched along the vertical direction toward the sea bottom using a stretching radio of 1.1.
  
  \begin{table}[htb]
  	\centering
  	\caption{The control parameters for the meshes used in the four different FEM methods that are implemented in this study to perform the hydrodynamic analyses.}
  	\begin{tabular}{ccccc}
  		\toprule
  		Method & $N_p$ & $N_{rx}$& $N_{ox}$ & $N_{oy}$\\
  		\midrule
  		Linear FEM     & 78526& 105& 300& 60    
  		  \\
  		Linear XFEM    & 81421& 105& 300& 60
  		  \\
  		Quad. FEM  & 15221  & 15& 120& 20
  		  \\
  		Quad. XFEM & 15416& 15& 120& 20 \\
  		\bottomrule
  	\end{tabular}
  	\label{tab3}
  \end{table}
  
  \par The added mass and damping coefficients are calculated by the four different FEMs, and the results are compared with the experimental results reported in \cite{1968The}, the linear numerical potential-flow calculations by \cite{liang2015application}. \cite{liang2015application} used the 2D HPC method in the frequency domain, and have taken account of the local singularity by a domain decomposition strategy, where the local corner-flow solutions were matched with the outer domain represented by the harmonic-polynomial cells. The mesh parameters used in our FEMs are listed in Table~\ref{tab3}, in which $N_p$ denotes the number of DOFs (including additional DOFs for XFEM) in the computational domain. The present numerical results agree excellently well with those by \cite{liang2015application}.  All numerical results seem to deviate from with the experimental results at low frequencies. As it is commented in \cite{1968The}, the uncertainties in the experimental results for $\omega^2B/(2g)<0.25$ may have been high. For $\omega^2B/(2g)\ge 0.25$, the numerical results agree better with the experiments. The small differences may have been contributed by the viscous flow separation at the sharp and other nonlinearities which will occur in reality. 
  
  From the results in Fig.~\ref{Fig.16}, we may conclude that all the numerical methods in the comparison are be able to accurately predict the linear hydrodynamic coefficients with an affordable effort. It is also observed that the XFEMs do not show clear advantages in the linear hydrodynamic analysis, which is expected as only integrals of velocity potential (multiplied by the normal vector) over the mean wetted body surface are involved in the pressure integration. As seen in the corner flow solution in Sect.~\ref{sect:corner-flow} , the velocity potential is not singular at the corner. However, the fluid velocity close to the sharp corners is singular, which poses great challenges in nonlinear wave loads analysis as will be explained further.  
  \begin{figure*}[t]
	\centering
	\subfigure[Linear XFEM]
	{\includegraphics[scale=.30]{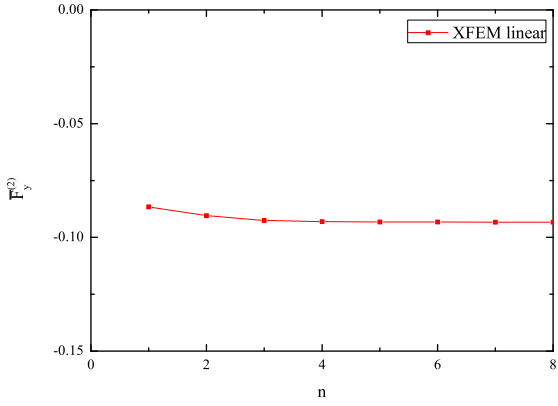}}
	\subfigure[Quadratic XFEM]
	{\includegraphics[scale=.30]{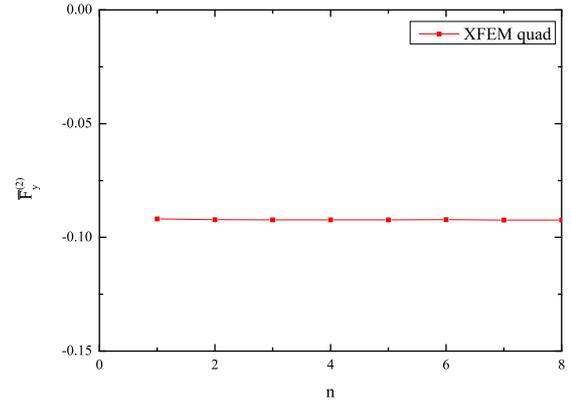}}
	\caption{Non-dimensional 2nd order mean vertical force versus the number of enrichment functions for non-dimensional oscillatory frequency of $\omega^2B/(2g)=1.0$. The non-dimensional 2nd order mean vertical force $\bar F_{y}^{(2)}=F_{y}^{(2)}/(\rho\omega^2\eta_{3a}^{2}B)$, $F_{y}^{(2)}$= 2nd order mean vertical force, $\rho$= mass density of water, $\eta_{3a}$= heave amplitude, $B$= beam, $n$= enrichment function number. Linear XFEM employs mesh 2 in Table~\ref{tab4}, quadratic XFEM employs mesh 2 in Table~\ref{tab5}.}
	\label{Fig.30}	
  \end{figure*}
  
  \begin{figure*}[t]
	\centering
	\subfigure[Linear XFEM]
	{\includegraphics[scale=.30]{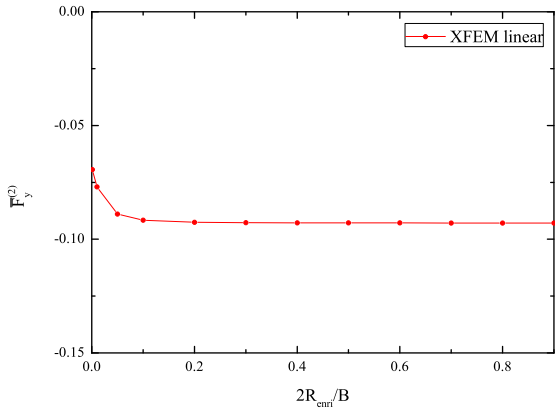}}
	\subfigure[Quadratic XFEM]
	{\includegraphics[scale=.30]{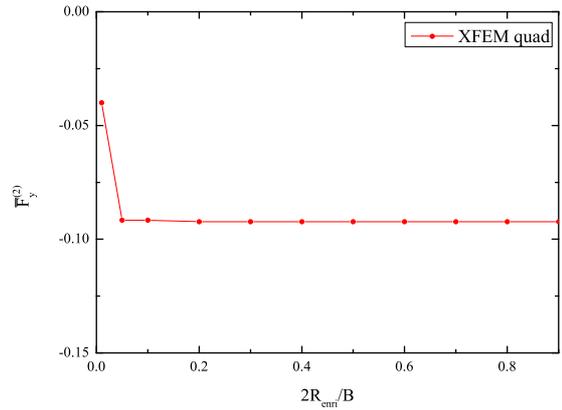}}
	\caption{Non-dimensional 2nd order mean vertical force versus enrichment radius. The considered non-dimensional oscillatory frequency is $\omega^2B/(2g)=1.0$. The non-dimensional 2nd order mean vertical force $\bar F_{y}^{(2)}=F_{y}^{(2)}/(\rho\omega^2\eta_{3a}^{2}B)$, $F_{y}^{(2)}$= 2nd order mean vertical force, $\rho$= mass density of water, $\eta_{3a}$= heave amplitude, $B$= beam, $r$= enrichment radius. Linear XFEM employs mesh 2 in Table~\ref{tab4}, quadratic XFEM employs mesh 2 in Table~\ref{tab5}.}
	\label{Fig.15}	
  \end{figure*}

  \begin{figure*}[ht]
  	\centering
  	\subfigure[Linear method]
  	{\includegraphics[scale=.30]{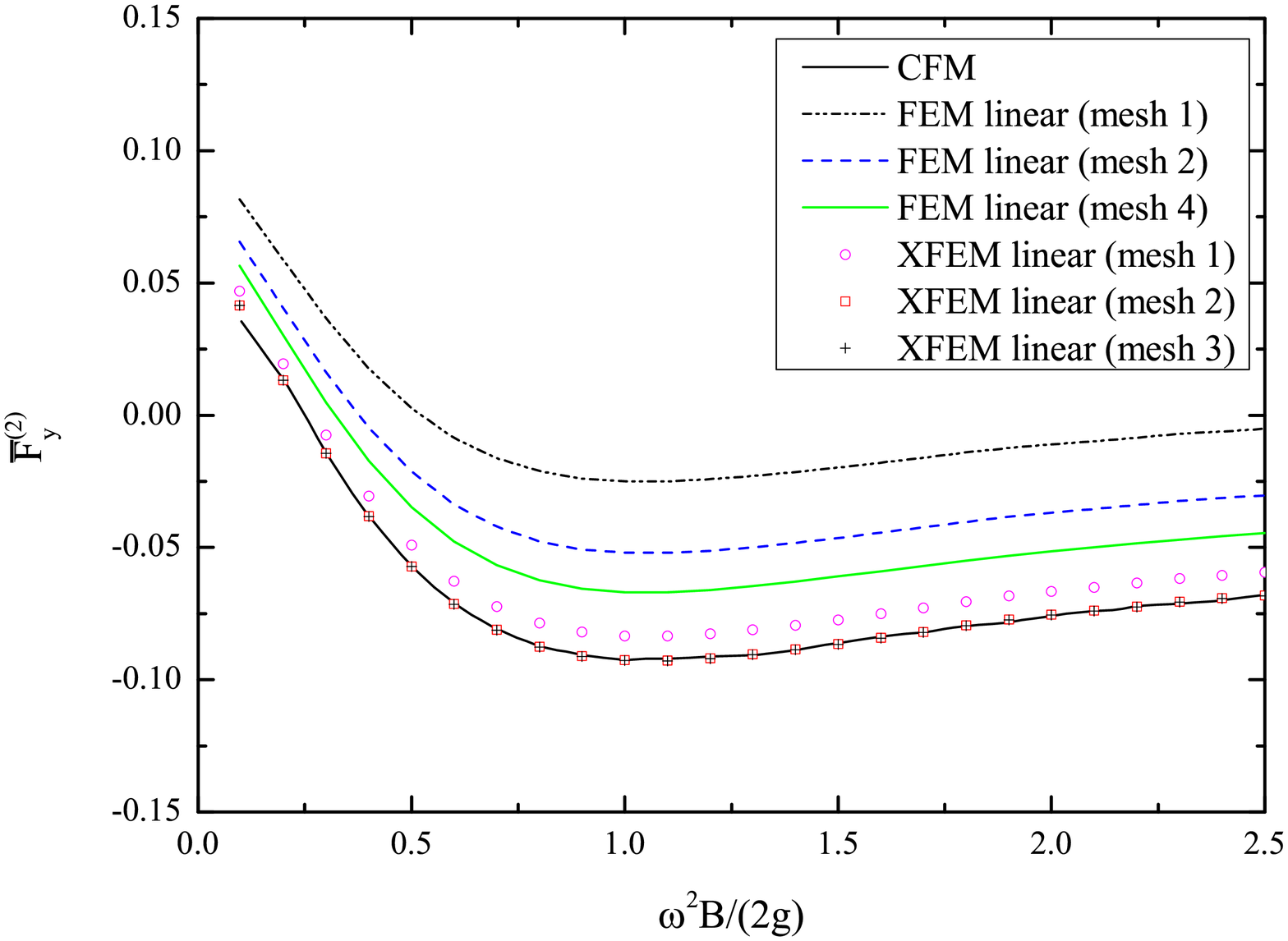}}
  	\subfigure[Quadratic method]
  	{\includegraphics[scale=.30]{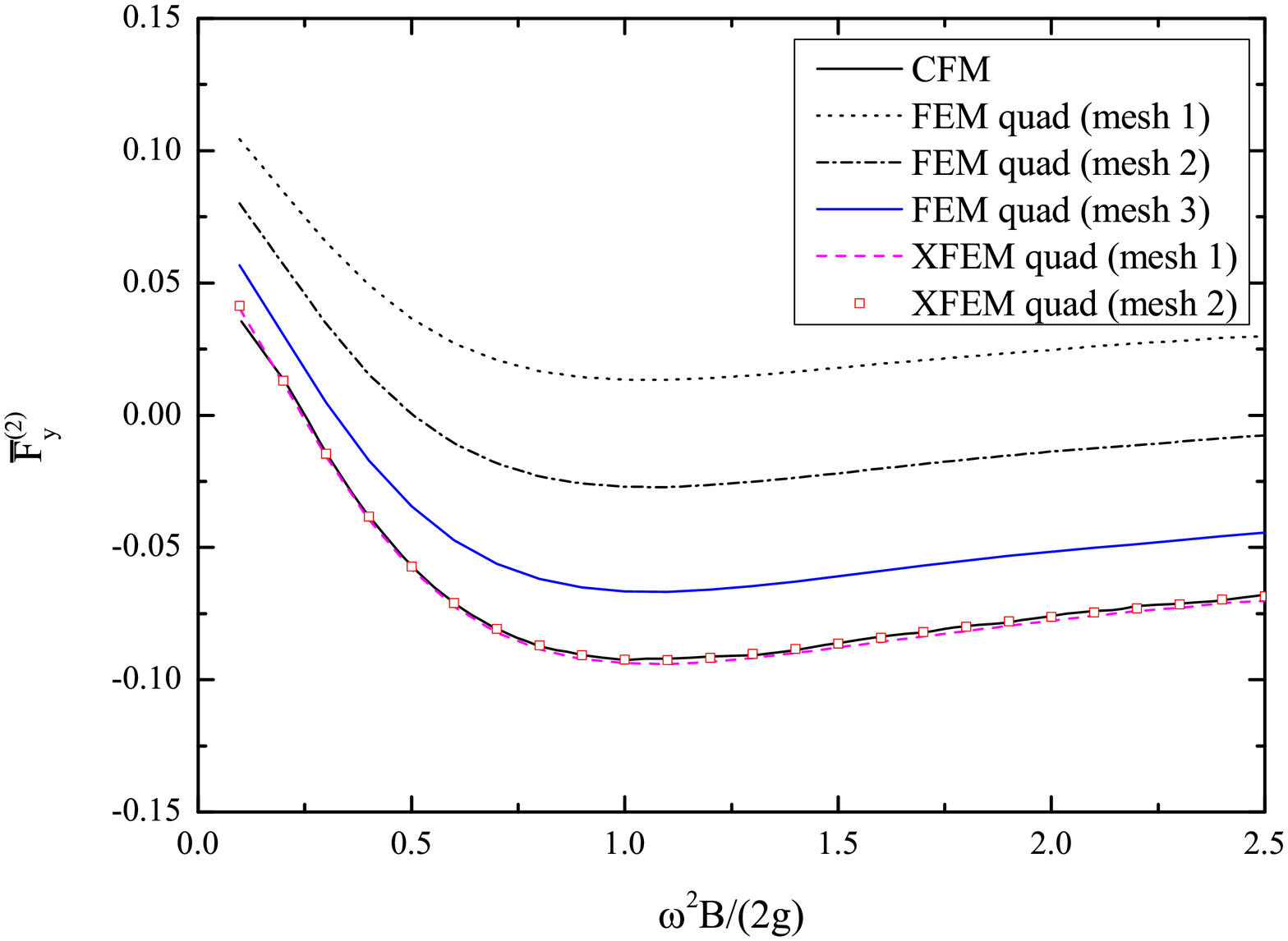}}
  	\caption{The non-dimensional 2nd order mean vertical force of a heaving floating rectangle. The non-dimensional 2nd order mean vertical force $\bar F_{y}^{(2)}=F_{y}^{(2)}/(\rho\omega^2\eta_{3a}^{2}B)$, $F_{y}^{(2)}$= 2nd order mean vertical force, $\rho$= mass density of water, $\eta_{3a}$= heave amplitude, $B$= beam, $\omega$= circular frequency. The CFM= conservation of fluid momentum. }
  	\label{Fig.17}	
  \end{figure*}
  
  \subsubsection{The 2nd order mean vertical force}
  
  \par The calculation of 2nd order wave loads based on pressure integration involves the integration of the quadratic terms of fluid velocities on body surface, which are singular but integrable near the sharp corners. \cite{zhao1989interaction} showed that the near-field approach based on direct pressure integration without special consideration of the singularity is very difficult to achieve convergent results, and the approach based on momentum and energy relationship or similar were much more efficient and robust. The later approach often involves integration on a control surface and a free surface confined by control surface and structure surface. Similar conclusions have been obtained later by many others \citep[e.g., see][]{chen2007middle, 2018Numerical, 2020A}.
  
  \par The time averaged 2nd order vertical hydrodynamic force acting on the heaving rectangle by direct pressure integration over the mean wet body surface can be expressed as: 
  \begin{equation}
  \label{Eq.42}
  \begin{aligned}
     \bar F_{y}^{\left( 2 \right)}  =& \frac{1}{T}\int_{0}^{T}\Bigg\{-\rho \int_{{{S}_{B0}}}\Bigg\{
          \frac{\partial {{\phi }^{\left( 1 \right)}}}{\partial t}
          + {{\eta }_{3}}\frac{{{\partial }^{2}}{{\phi }^{\left( 1 \right)}}}{\partial y\partial t} +\\
         &  \frac{\partial {{\phi }^{\left( 2 \right)}}}{\partial t}
          + \frac{1}{2}{{\left[ \frac{\partial {{\phi }^{\left( 1 \right)}}}{\partial x} \right]}^{2}
          + \frac{1}{2}{{\left[ \frac{\partial {{\phi }^{\left( 1 \right)}}}{\partial y} \right]}^{2}} } \Bigg\} {{n}_{3}}\text{d}S \Bigg\}\text{d}t,
  \end{aligned}
  \end{equation}
  where $S_{B0}$ denotes the wetted mean body surface. $\eta_{3}$ is the heave motion define as $\eta_{3} =\Re[\eta_{3a}\ee^{\ii\omega t}]$, with $\eta_{3a}$ as the heaving amplitude. $n_3$ is the vertical component of the normal vector. $T$ is the oscillatory period expressed as $T = 2\pi/\omega$. $\phi^{(1)}$ and $\phi^{(2)}$ represent the first and second order velocity potential, respectively. A waterline integral due the fluctuation of waves near the mean water level is neglected as it does not contribute to the vertical loads in this particular case. The time derivatives of first and second order velocity potential equal to zero after time average over one period, and thus Eq.(\ref{Eq.42}) can be simplified as:
  \begin{equation}
  \label{43}
    \bar F_{y}^{(2)} =  -\rho\overline{\int_{{{S}_{B0}}} \left[ {{\eta}_{3}}\frac{{{\partial }^{2}}{{\phi }^{( 1)}}}{\partial y\partial t} + 
    \frac{1}{2}\nabla\phi^{(1)}\cdot\nabla\phi^{(1)} \right] {{n}_{3}}\text{d}S}.
  \end{equation}
  
  \par From a theoretical perspective, the near-field approach, far-field approach and the approaches based on control surfaces should be mathematically equivalent. However, since it is very difficult for the conventional numerical methods, e.g. FEM, FDM and BEM, to accurately describe the exact fluid velocities close to sharp corners, slow grid-convergences are expected for the near-field approach when it is applied to calculate the 2nd order wave loads. Despite difficult, it is still believed by the authors of this paper that, the strong variation of the local velocities can be captured accurately if an appropriate numerical method is adopted, thus the near-field approach can still be a good option for 2nd wave-load analysis. A good example of such a numerical method is the domain decomposition strategy developed by \cite{liang2015application}, where the solutions in the local domain surrounding the sharp corners are represented by the analytical corner-flow solutions. The strategy leads to very accurate and efficient near-field result, but it is not easy to implement for general purposes. The XFEM is a more powerful and general-purpose framework, which allows us to easily and explicitly include, for instance the singular corner-flow solutions as enrichment functions, in the local finite-element approximations. It also inherits other good features of the conventional FEMs, e.g. the use of unstructured meshes. 
  
  \par For the considered rectangle, the interior angle at each corner is $\gamma = 90 ^{\circ}$, where $\gamma$ is the interior angle as illustrated in Fig.~\ref{Fig.1}.  Eq.~\eqref{Eq.9} presents all possible fundamental solutions to the corner flows, among which we choose only the first a few as our enrichment function. The first term with $j=1$ is $\varphi = A_1 r^{\frac{2}{3}} \cos (\frac{2}{3}\theta)$, and the resulting radial velocity $\frac{\partial \varphi}{\partial r}$ and circumferential velocity $\frac{1}{r} \frac{\partial \varphi}{\partial \theta} $ are in the form of $r^{-\frac{1}{3}}-$singularity as $r\rightarrow 0$, which are difficult for any regular functions to achieve good approximation. 
  
  In Fig.~\ref{Fig.30}, we compare the non-dimensional $\bar F_{y}^{(2)}$ when different number of terms from Eq.~\eqref{Eq.9} are included as enrichment functions. 
   As shown in Fig.~\ref{Fig.30}(a), for linear XFEM, convergent result has been achieved for  enrichment function number $n\ge 3$. For quadratic XFEM, the convergence will be achieved with $n\ge 1$, as demonstrated in Fig.~\ref{Fig.30}(b). The reason that a linear XFEM needs more enrichment functions than a quadratic XFEM is as follows: the fundamental solution of a corner flow contains a singular term with $j=1$ in Eq.~\eqref{Eq.9} and other higher-order non-singular terms with $j\ge 2$. Those non-singular terms are more accurately captured by the regular quadratic shape functions, and thus it seems to be sufficient for a quadratic XFEM to use only the singular enrichment function from Eq.~\eqref{Eq.9}. Based on the discussion above and the numerical observation, only three enrichment functions will be considered in later analyses. Adding unnecessarily too many higher-order terms with $j>3$ will pose extra difficulties in numerical integration. On the other hand, the extra DOFs due to enrichment will increase rapidly with the number of enrichment function at each enrichment point.
  
  \par Fig.~\ref{Fig.15} displays the non-dimensional 2nd order mean vertical force $\bar F_{y}^{(2)}$ for $\omega^2B/(2g)=1.0$ as function of $2R_{enri}/B$. The numerical results indicate that, for both linear and quadratic XFEMs, the convergence is achieved when $2R_{enri}/B$ $\geq 0.2$. The results also suggest that it is unnecessary to use a too large enrichment radius, because the results do not seem to improve further as long as $R_{enri}$ is greater than a threshold value of approximately $0.2$. On the other hand, larger $R_{enri}$ also means more extra DOFs and unknowns. 
  
  \par In Fig.~\ref{Fig.17}(a), the numerical results of $\bar F_{y}^{(2)}$ by the linear FEM and the linear XFEM are compared with a reference solution in \cite{liang2015application} based on conservation of fluid momentum (CFM). Direct pressure integration has been applied in the present FEM analyses. Mesh 1, mesh 2 and mesh 3 in the parentheses indicate coarse, medium and fine meshes, respectively. Details of the mesh parameters are shown in Table~\ref{tab4}. Apparently, the linear XFEM is more accurate than linear FEM as seen from their comparisons with the CFM results \citep{liang2015application}. Convergent result can be achieved rapidly after refine mesh 1 to mesh 2 via linear XFEM. The unknown numbers (or total DOFs) of mesh 1 and mesh 2 are 28866 and 81421 respectively when linear XFEM is applied. On the contrary, the linear FEM has not reached the convergence even with the finest mesh, i.e. mesh 4 with total DOFs of $N_p=556146$ in Table~\ref{tab4}. Note that the total number of unknowns, or total DOFs, are different for a FEM and a XFEM, even though the same mesh is used. This is due to the extra DOFs introduced in the XFEM as a result of local enrichment.
  \begin{table}[htb]
  	\centering
  	\caption{The three different meshes and DOFs parameters for the two linear (FEM and XFEM) methods, which are used in the calculation of the 2nd order mean vertical force.}
  	\begin{tabular}{cccccc}
  		\toprule
  		Method & $N_p$ & $N_{rx}$& $N_{ox}$ & $N_{oy}$\\
  		\midrule
  		Linear FEM (mesh 1)& 28275& 25& 300& 60    
  		  \\
  		Linear FEM (mesh 2)& 78526& 105& 300& 60
  		  \\
  		Linear FEM (mesh 4)& 556146& 405& 400& 80
  		  \\
  		Linear XFEM (mesh 1)& 28866& 25& 300& 60
  		  \\
  		Linear XFEM (mesh 2)& 81421& 105& 300& 60 
  		\\
  		Linear XFEM (mesh 3)& 99084& 125& 300& 60 
  		\\
  		\bottomrule
  	\end{tabular}
  	\label{tab4}
  \end{table}

  \par For quadratic methods, we also consider three different meshes, i.e. coarse, medium and fine meshes, represented by mesh 1, 2 and 3 in Table~\ref{tab5}, respectively. As illustrated by the comparisons in Fig.\ref{Fig.17}(b), the conventional quadratic FEM is not able to reach a convergence even with the finest mesh (mesh 3) with a total DOFs of $N_p=406631$. On the contrary, the quadratic XFEM results are convergent with the medium mesh (mesh 2, $N_p=15416$). In fact, results based on the coarse mesh (mesh 1, $N_p=9293$) are already very close to the reference results. In this coarse mesh resolution, only 4 elements are distributed on half of the rectangle bottom. 
  
  \begin{table}[htb]
  	\centering
  	\caption{Mesh parameters for the two quadratic (FEM and XFEM) methods, which are applied in the calculation of the 2nd order mean vertical force.}
  	\begin{tabular}{ccccc}
  		\toprule
  		Method & $N_p$ & $N_{rx}$& $N_{ox}$ & $N_{oy}$\\
  		\midrule
  		Quad. FEM (mesh 1) & 9281& 4& 120& 20    
  		  \\
  		Quad. FEM (mesh 2) & 15221& 15& 120& 20 
  		  \\
  		Quad. FEM (mesh 3)& 406631& 215& 120& 50 
  		  \\
  		Quad. XFEM (mesh 1) & 9293& 4& 120& 20 
  		  \\
  		Quad. XFEM (mesh 2)& 15416& 15& 120& 20  \\
  		\bottomrule
  	\end{tabular}
  	\label{tab5}
  \end{table}

  \par Comparing the two XFEMs, quadratic XFEM has shown much faster mesh-convergence rate than linear XFEM. More specifically, convergent results can be reached by quadratic XFEM with less than $N_p=15416$ DOFs, while its takes $N_p=81421$ for the linear XFEM. Therefore, the quadratic XFEM is considered as more competitive. From the standpoint of solution enrichment, the quadratic XFEM can be seen as a combination of global and local enrichment, with a global enrichment achieved via higher Lagrange polynomials in regular shape functions, and a local enrichment realized by adding prior knowledge to the local approximation space. The linear XFEM, however, only enriches the solution locally. Therefore, it is generally expected that the quadratic XFEM over-performs the linear XFEM.
  
  \begin{figure}[t]
		\centering
		\includegraphics[scale=.3]{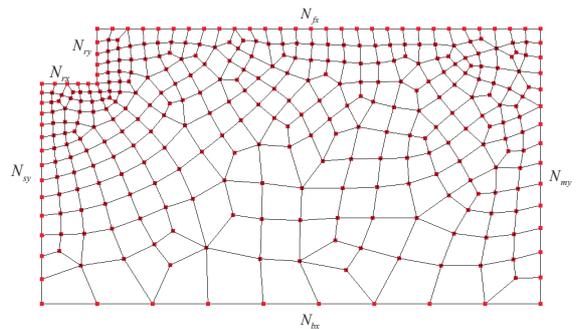}
		\caption{An example of the unstructured mesh of linear elements for the half rectangular heaving on the free surface. $N$ with subscripts represent the number of element along the boundaries of the fluid domain.}
		\label{Fig.31}	
	\end{figure}
	
	\subsubsection{Application of unstructured meshes}
    \begin{figure*}[ht]
    	\centering
    	\subfigure[Linear method]
    	{\includegraphics[scale=.30]{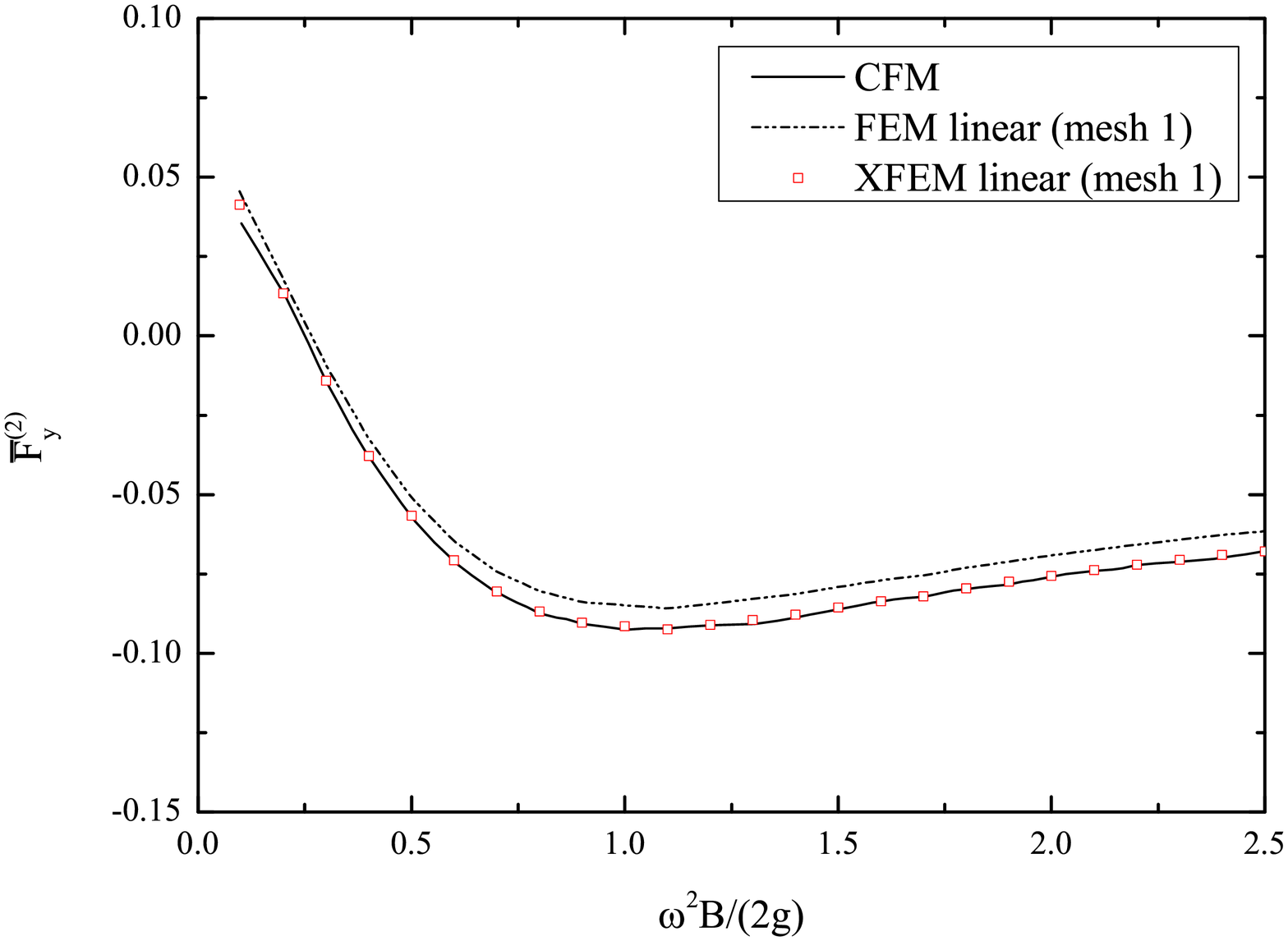}}
    	\subfigure[Quadratic method]
    	{\includegraphics[scale=.30]{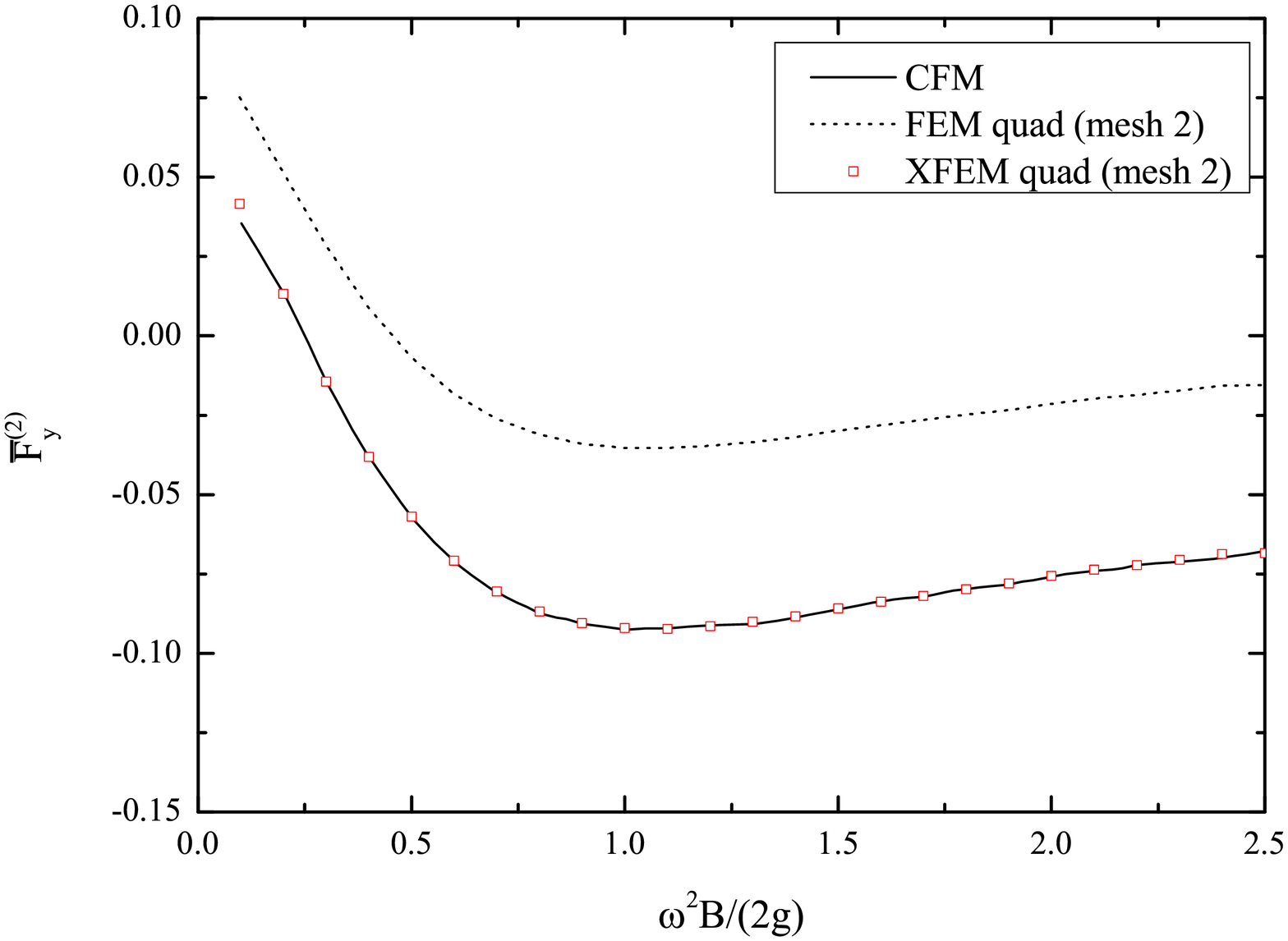}}
    	\caption{The non-dimensional 2nd order mean vertical force of a heaving floating rectangle with unstructured mesh. The non-dimensional 2nd order mean vertical force $\bar F_{y}^{(2)}=F_{y}^{(2)}/(\rho\omega^2\eta_{3a}^{2}B)$, $F_{y}^{(2)}$= 2nd order mean vertical force, $\rho$= mass density of water, $\eta_{3a}$= heave amplitude, $B$= beam, $\omega$= circular frequency. The CFM= conservation of fluid momentum. }
    	\label{Fig.32}	
    \end{figure*}
    
    \begin{table*}[htb]
    	\centering
    	\caption{Mesh parameters for the conventional linear and quadratic FEMs and their corresponding XFEMs, which are applied to obtain the 2nd order mean vertical force.}
    	\begin{tabular}{cccccccc}
    		\toprule
    		Method & $N_p$ & $N_{rx}$& $N_{ry}$& $N_{fx}$ & $N_{bx}$ &  $N_{sy}$ &  $N_{my}$\\
    		\midrule
    		Linear FEM (mesh 1) & 35565& 75& 75& 1199&   69&  59& 59
    		\\
    		Linear XFEM (mesh 1) & 40854& 75& 75& 1199&   69&  59& 59
    		\\
    		Quad. FEM (mesh 2) & 5756& 10& 10& 199 & 29& 19& 19
    		\\
    		Quad. XFEM (mesh 2) & 5870& 10& 10& 199 & 29& 19& 19\\
    		\bottomrule
    	\end{tabular}
    	\label{tab6}
    \end{table*}

	\par In the previous subsections, a multi-block structured mesh was adopted for demonstration purpose, and the numerical results based on XFEMs were very encouraging. However, it is well-known that, one of the most powerful property of FEM is that it allows for the use of unstructured mesh without having to modify the numerical code. It is much easier for the unstructured meshes to deal with problems involving complex boundaries. In this subsection, the unstructured mesh will be adopted to investigate the same problem that have been studied in the previous subsection. 
	
	\par An example of the unstructured mesh close the 2D rectangle, generated from the open-source mesh generator GMSH, is shown in Fig.~\ref{Fig.31}. The following parameters are defined to control the number of elements on the fluid boundaries: $N_{rx}$ is the number of elements on the bottom of the rectangle, $N_{ry}$ along the side wall of the rectangle, $N_{fx}$ along the free surface, $N_{sy}$ along the symmetry face, $N_{bx}$ along the bottom of the computational domain and $N_{my}$ along the matching boundary. Furthermore, for both linear and quadratic mesh, the mesh is stretched by a fixed stretching radio of 1.1 along the body boundary, so that the meshes are finer close to the corners. The meshes are also stretched vertically towards the bottom of the fluid and horizontally towards the matching boundary, using stretching factors of 1.08 and 1.05 respectively. The meshes are so adapted that the mesh density is higher close to the body and the free surface. 
	
	\par The 2nd order mean vertical force on the heaving rectangle at free surface is studied again in the frequency domain by using the unstructured mesh and the four FEMs, and the corresponding results for linear FEMs and quadratic FEMs are shown in Fig.~\ref{Fig.32}(a) and Fig.~\ref{Fig.32}(b) respectively. The main parameters of the applied unstructured meshes are summarized in Table~\ref{tab6}. 
	
	Due to the use of unstructured meshes and stretched grid on the fluid boundaries, it is expected that the required total number of unknowns are much smaller than that of the multi-block structured meshes. This has also been confirmed by our numerical results in Fig.~\ref{Fig.32}(a) and Fig.~\ref{Fig.32}(b). As seen in the figures, to achieve convergent results for $\bar F_{y}^{(2)}$, it is sufficient to use mesh 1 (total DOFs $N_p=40854$) and mesh 2 ($N_p=5870$) in Table~\ref{tab6} for linear XFEM and quadratic XFEM, respectively. On the other hand, as expected, the results of the conventional FEMs are not convergent when the same meshes as the corresponding XFEMs are used. There is one point that must be clarified for the results of the conventional linear FEM and quadratic FEM. In Fig.~\ref{Fig.32}, the linear FEM results appear to be closer to reference results than that of quadratic FEM. This is due to the fact that, mesh 1 as used by linear FEM is much finer than mesh 2 used by quadratic FEM. 
	
	\begin{figure*}[t]
  \centering
	\subfigure[Linear element] {\includegraphics[scale=.80]{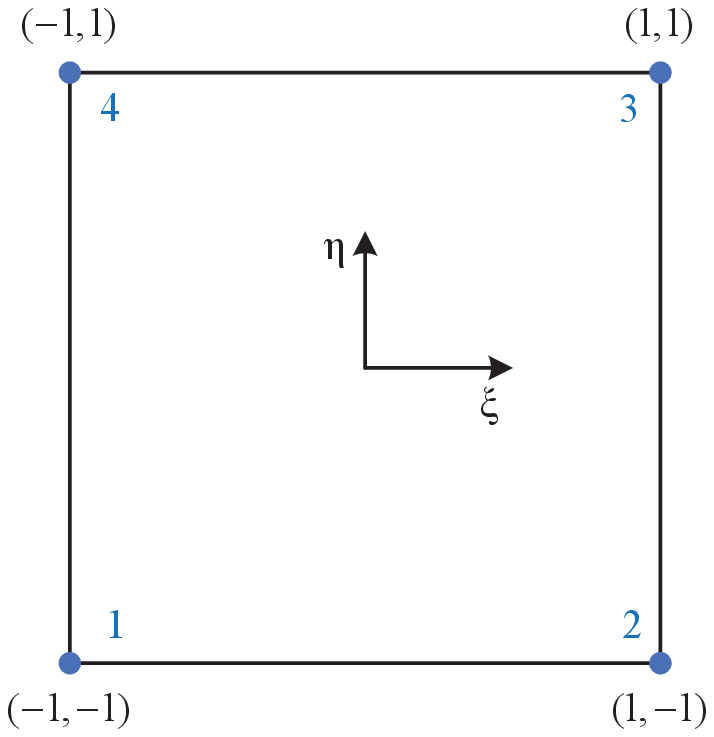}}
	\subfigure[Quadratic element] {\includegraphics[scale=.80]{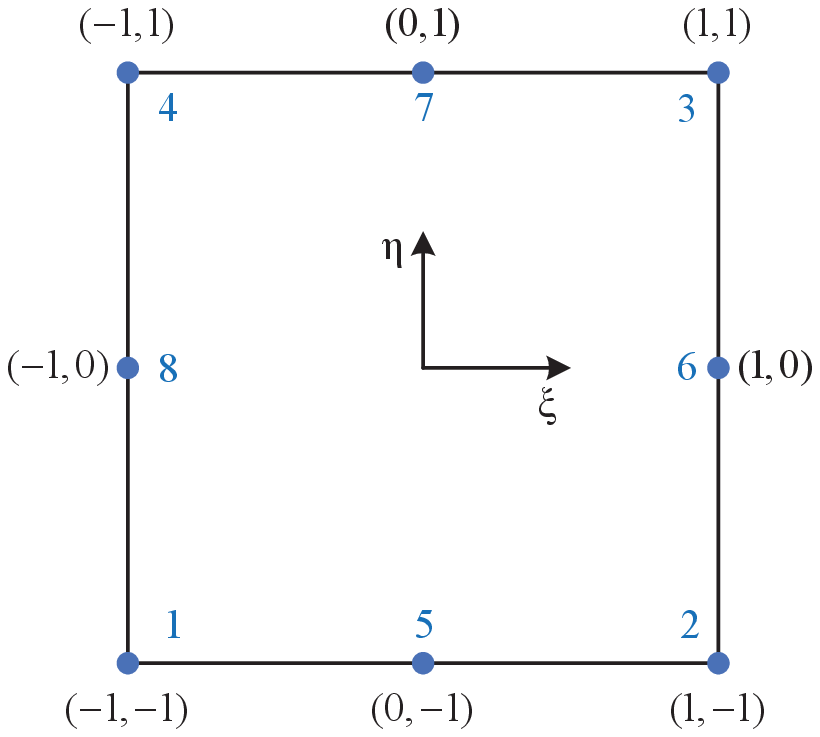}}
	\caption{Linear and quadrilateral standard element on a $\xi\eta$- plane.}
	\label{Fig.22}	
 \end{figure*}
 \begin{figure}[t]
  	\centering
  	\includegraphics[scale=.4]{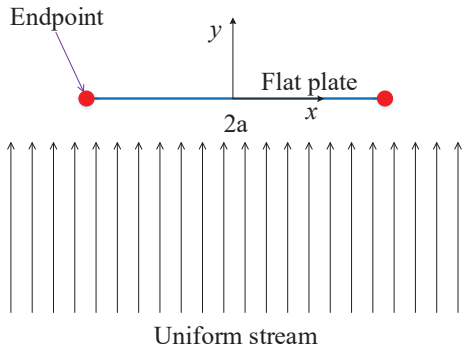}
  	\caption{The uniform flow around the flat plate.}
  	\label{Fig.18}	
  \end{figure}
  \begin{figure*}[t]
  	\centering
  	\subfigure[Without modification] {\includegraphics[scale=.50]{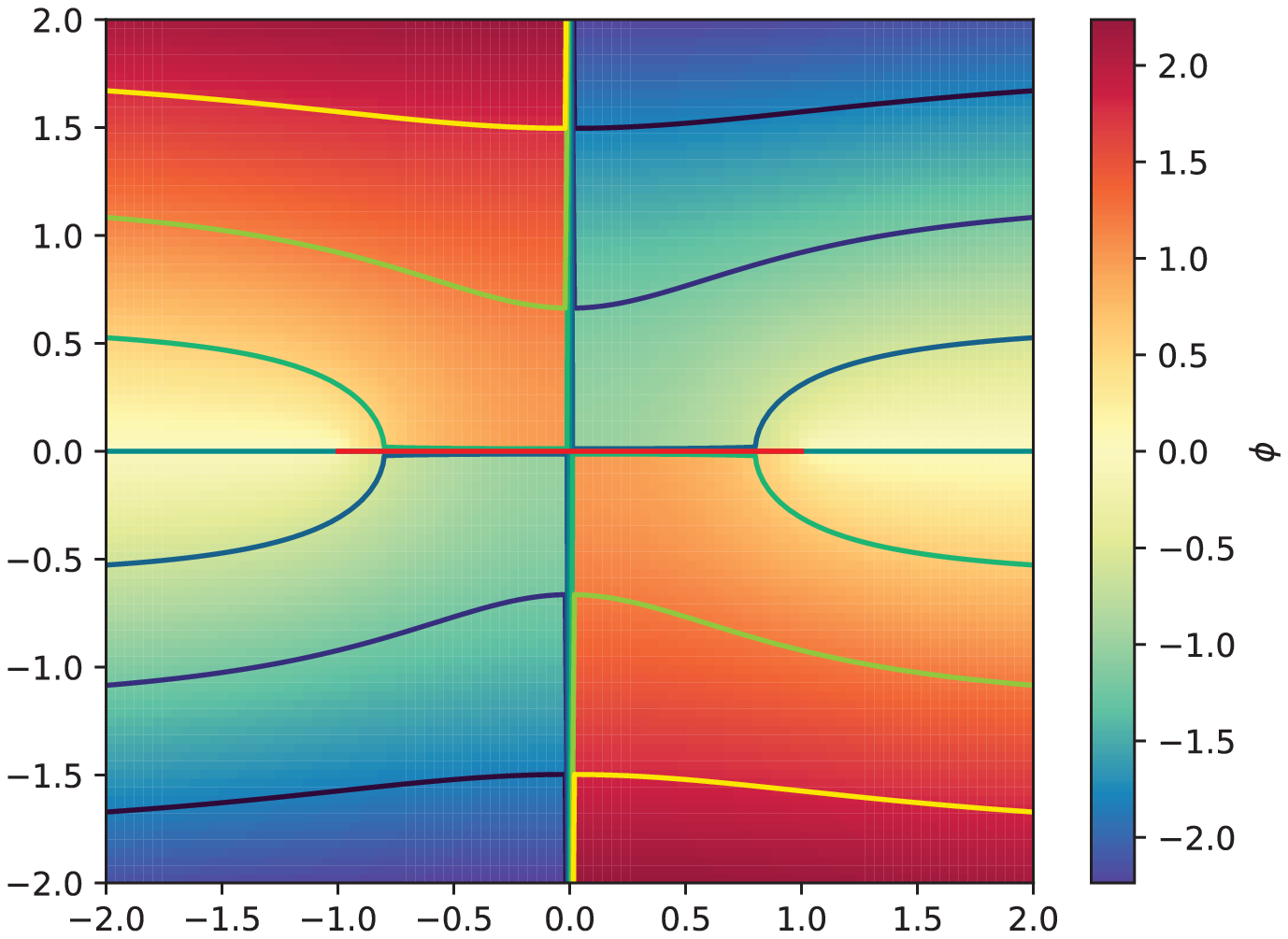}}
  	\subfigure[Modification] 
  	{\includegraphics[scale=.50]{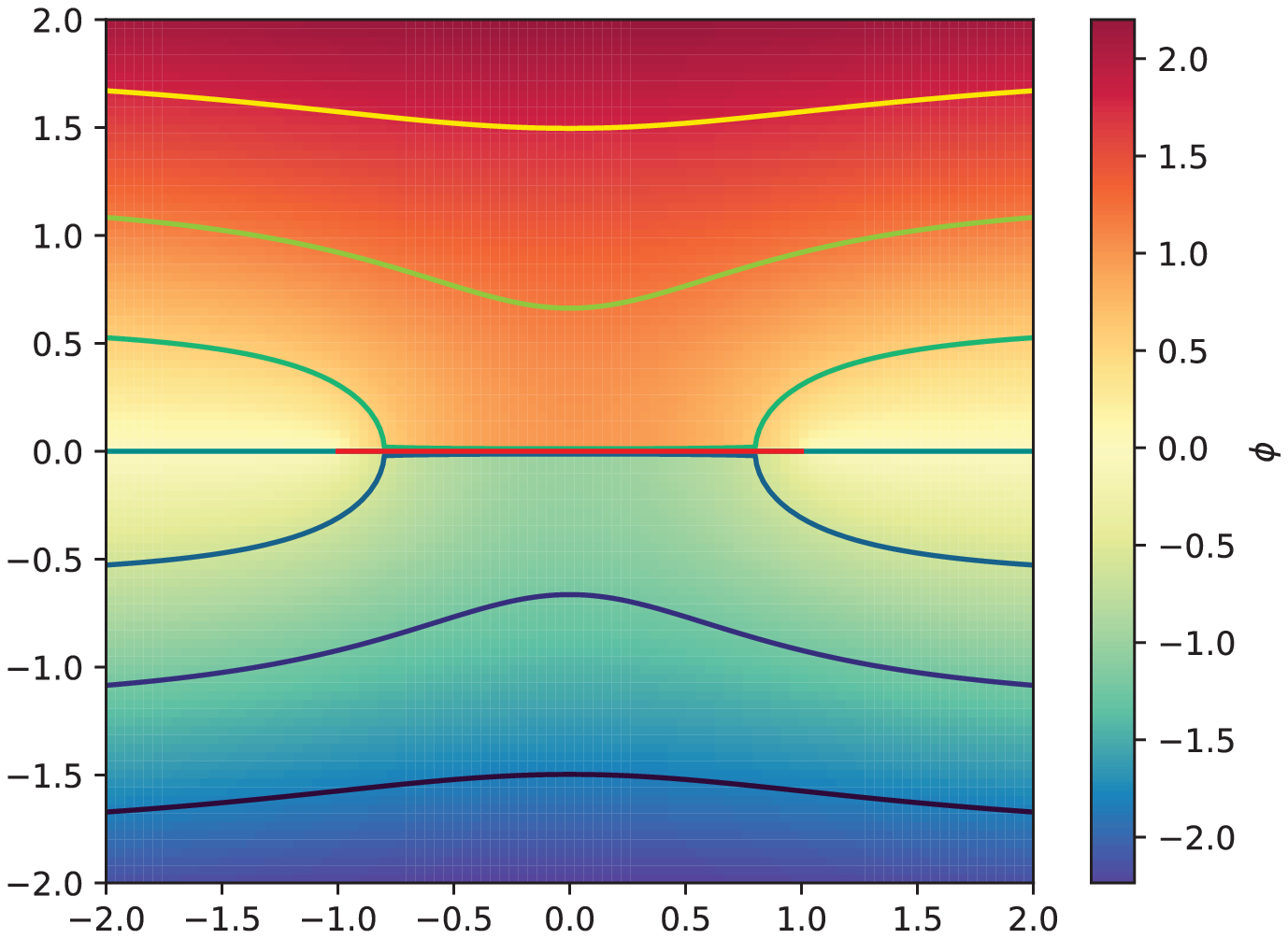}}
  	\caption{The contour of velocity potential.}
  	\label{Fig.19}	
  \end{figure*}
  \begin{figure*}[t]
  	\centering
  	\subfigure[Without modification] {\includegraphics[scale=.50]{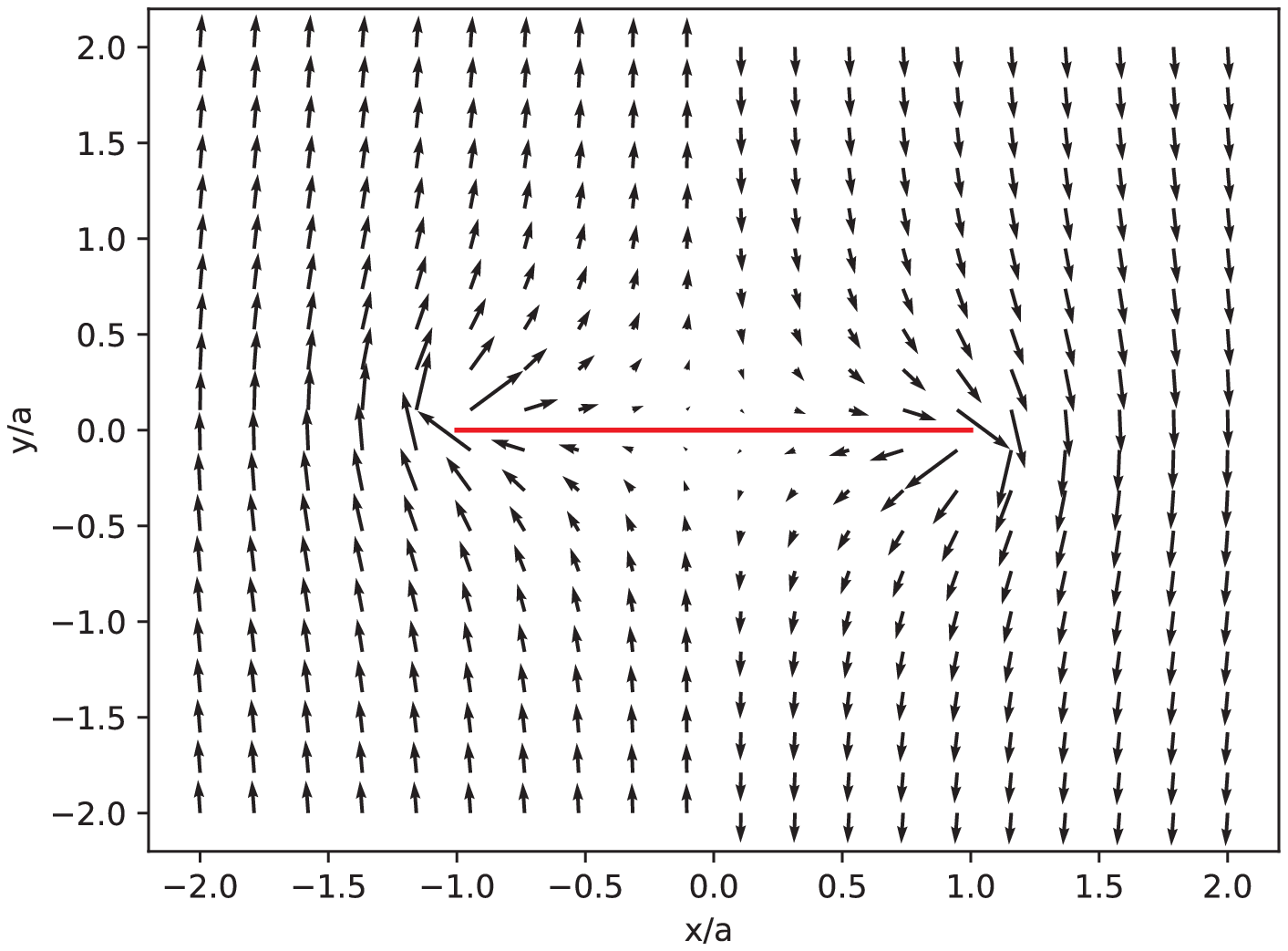}}
  	\subfigure[Modification] 
  	{\includegraphics[scale=.50]{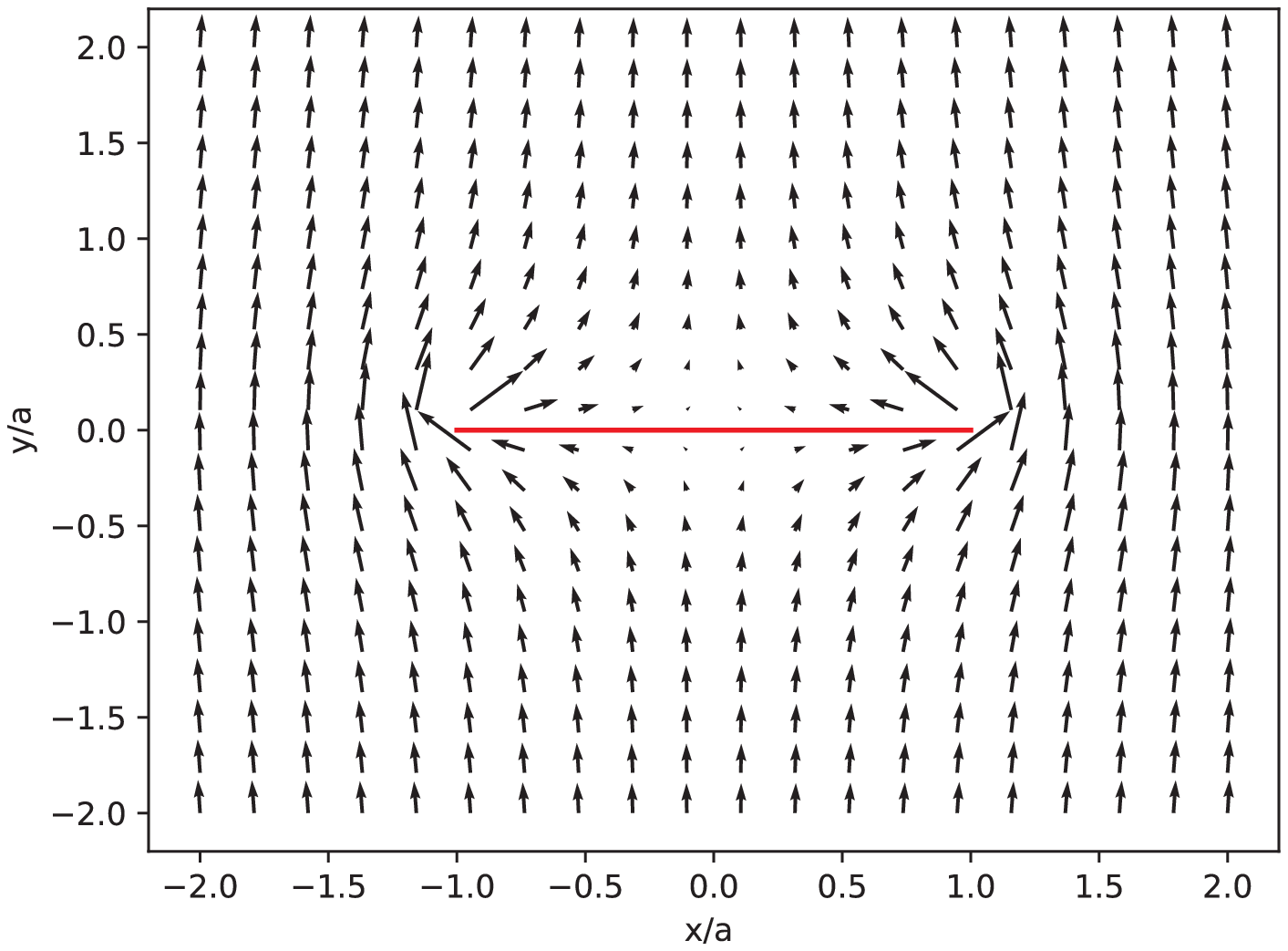}}
  	\caption{Velocity vector diagram.}
  	\label{Fig.20}	
  \end{figure*}
  
	In \cite{liang2015application}, a modified HPC method based on domain decomposition method was developed to solve the same hydrodynamic problem of the heaving rectangle at free surface in the frequency domain. Corner-flow solutions were used in the inner domain surrounding the sharp corner, while the outer domain solutions were represented by overlapping harmonic-polynomial cells. The inner and outer domain solutions are matched at their common boundaries. This method was shown to be capable of providing convergent 2nd order mean wave loads by using 80 elements along half of the bottom and in total approximately 352000 unknowns, while our linear and quadratic XFEM models need much fewer unknowns (only 5870 for quadratic XFEM and 40854 for linear XFEM) to achieve equally good results. In a nutshell, the superiority of the present study over \cite{liang2015application} is twofold. Firstly, from implementation point of view, the enrichment strategy based on Partition of Unity in XFEMs to include singular functions near the corner is easier and more flexible. Secondly, the unstructured meshes are allowed in XFEMs, which enables XFEMs to deal with more complex structures, whereas it is expected to be more difficult for the HPC method.

\section{Conclusions}\label{sect:Conclusion and perspective}
  \par The XFEM is applied as an accurate and efficient tool to solve 2D potential-flow hydrodynamic problems for structures with sharp edges. To demonstrate the advantages of XFEM, four FEM codes, in 2D, including the conventional linear and quadratic FEMs, and the two corresponding XFEMs, are implemented and compared. All of our results have confirmed that the XFEM is a promising framework to deal with potential-flow hydrodynamic problems involving structures with sharp edges. Three different enrichment strategies, including: the point enrichment, patch enrichment and radius enrichment, are also investigated in the study of the uniform flow around an infinite-thin flat plate. The first two enrichment methods are found to be mesh-dependent, and are not able to achieve the expected spatial convergence rate. 
  However, the radius enrichment method is mesh-independent, and has shown remarkably better accuracy and spatial convergence rate. Therefore, it is considered as the best option over the other three counterparts. By studying the horizontal fluid velocity along the flat plate, we also demonstrate that XFEMs are capable of capturing the strong flow variation close to the endpoints, which cannot be represented by the conventional FEMs.
  
  \par For a heaving rectangular cylinder on the free surface, both the conventional FEMs and XFEMs can accurately predict the linear hydrodynamic coefficients with acceptable computational efforts, indicating that the singularity at sharp corner is inconsequential to the linear hydrodynamic loads. However, it has important effects on the 2nd order mean wave loads if the direct pressure integration is employed, because the singular flow velocities are involved. Compared with reference results based on conservation of fluid momentum, both linear and quadratic XFEMs have shown encouraging results even with a relative coarse mesh resolution, while the quadratic XFEM has an overall better performance than the linear XFEM. On the contrary, it is difficult for the two conventional FEMs to achieve convergence even with an extremely fine mesh.
  
  For the quadratic XFEM, it is also found sufficient to only include the first singular term from the corner-flow solutions in the local enrichment, while it is beneficial to include a few more, e.g. 3 terms, in the local enrichment in the linear XFEM. 
  
  As a final demonstration, we show that the adoption of unstructured meshes and local refinement close to the sharp edges have a great potential to further reduce the total number of unknowns to achieve a desired accuracy. 
  
\section*{Appendix A. shape function and isoparametric element}
  \par Commonly, the shape function is defined in an element, for simplicity, using $n^e_p$ to represent the number of the nodes in a single element. Referring to \cite{zienkiewicz2005finite}, for 4-node quadrilateral linear FEM, namely $n^e_p=4$, the shape function defined on a parametric $\xi\eta$-plane can be written as:
  \begin{equation}
	\label{Eq.54}
	  {{N}_{i}}=\frac{1}{4}\left( 1+{{\xi }_{i}}\xi  \right)\left( 1+{{\eta }_{i}}\eta  \right), \qquad i=1,\cdots,4.
  \end{equation}
  where $(\xi_i,\eta_i)$ denote the normalized coordinates at node $i$.
  For an incomplete quadratic quadrilateral element, namely $n^e_p=8$, the shape function can be expressed as:
  \begin{equation}
  \label{Eq.55}
  \begin{aligned}
    {{N}_{i}}& =\frac{1}{4}\left( 1+{{\xi }_{i}}\xi  \right)\left( 1+{{\eta }_{i}}\eta  \right)\left( {{\xi }_{i}}\xi +{{\eta }_{i}}\eta -1 \right) \quad (i=1,\cdots,4),\\
	{{N}_{i}}& =\frac{1}{4}\left( 1+{{\xi }^{2}} \right)\left( 1+{{\eta }_{i}}\eta  \right)\quad(i=5,7),\\
	{{N}_{i}}& =\frac{1}{2}\left( 1+{{\xi }_{i}}\xi  \right)\left( 1-{{\eta }^{2}} \right)\quad(i=6,8).
  \end{aligned}
  \end{equation}
  Examples of the 4-node and 8-node quadrilateral elements in the parametric $\xi\eta$-plane are illustrated in Fig.~\ref{Fig.22}.
  
\section*{Appendix B. Analytical solution of flow over a flat plate}
  The complex potential of uniform around flat plate in a 2D infinite domain is \citep{Newman2017Marine}
  \begin{equation}
  \label{Eq.44}
    W\left( z \right)=-z{{v}_{0}}\cos \alpha +\ii{{v}_{0}}\sqrt{{{z}^{2}}-{{a}^{2}}}\sin \alpha,
  \end{equation}
  where $z=x+\ii y$,
  and the corresponding complex velocity is
  \begin{equation}
  \label{Eq.45}
    u-\ii v=-{{v}_{0}}\cos \alpha +\ii {{v}_{0}}\frac{z}{\sqrt{{{z}^{2}}-{{a}^{2}}}}\sin \alpha,
  \end{equation}
  where $v_0$ denotes the velocity of the uniform stream, $\alpha$ the angle between the uniform flow and the plate, $a$ the half of width of the flat plate, and $u$ and $v$ the horizontal and vertical velocity components, respectively. For convenience, we let $\alpha=\pi/2$ and $v_0=1$ as shown in Fig.~\ref{Fig.18}, the complex potential in the fluid domain is simplified to
  \begin{equation}
  \label{Eq.46}
    W(z)=\ii\sqrt{{{z}^{2}}-{{a}^{2}}},
  \end{equation}
  and the complex velocity become
  \begin{equation}
  \label{Eq.47}
    u-\ii v=\ii\frac{z}{\sqrt{{{z}^{2}}-{{a}^{2}}}}.
  \end{equation}
  The complex potential can be divided into two parts, including: the potential function $\phi(x,y)$, and the stream function $\chi(x,y)$
  \begin{equation}
  \label{Eq.48}
    W(z)=\phi( x,y)+\ii\chi( x,y).
  \end{equation}
  According to Eqs.~\eqref{Eq.46}) and \eqref{Eq.48}, we obtain
  \begin{equation}
  \label{Eq.49}
    \phi =\operatorname{Re}\left\{ W(z) \right\}=\operatorname{Re}\left\{ -\ii\sqrt{{{z}^{2}}-{{a}^{2}}} \right\}.
  \end{equation}
  The velocity in the fluid domain can be written as:
  \begin{equation}
  \label{Eq.50}
    u=\frac{\partial \phi }{\partial x} =\Re\left(\ii\frac{z}{\sqrt{{{z}^{2}}-{{a}^{2}}}}\right).
  \end{equation}
  \begin{equation}
  \label{Eq.51}
    v=\frac{\partial \phi }{\partial y}=-\Im\left(\ii\frac{z}{\sqrt{{{z}^{2}}-{{a}^{2}}}}\right)
  \end{equation}
  The velocity potential and velocity determined by Eqs.~\eqref{Eq.49} and \eqref{Eq.51} are not physical on the right-half plane as shown in Fig.~\ref{Fig.19} (a) and Fig.~\ref{Fig.20} (a), respectively. The velocity potential $\phi$ and velocity vectors are not symmetric about y-axis and there is discontinuity along $y$-axis.  The correction, therefore, should be made when $x\ge 0$ gives
  \begin{equation}
  \label{Eq.52}
  \begin{aligned}
    \phi =-\phi,  \\ 
     u=-u, \\ 
     v=-v.
  \end{aligned}
  \end{equation}
  After modification, the contour of velocity potential was shown in Fig.~\ref{Fig.19} (b) and the vector diagram of velocity in Fig.~\ref{Fig.20} (b).
  
\bibliography{References}

\begin{thebibliography}{53}
\expandafter\ifx\csname natexlab\endcsname\relax\def\natexlab#1{#1}\fi
\providecommand{\url}[1]{\texttt{#1}}
\providecommand{\href}[2]{#2}
\providecommand{\path}[1]{#1}
\providecommand{\DOIprefix}{doi:}
\providecommand{\ArXivprefix}{arXiv:}
\providecommand{\URLprefix}{URL: }
\providecommand{\Pubmedprefix}{pmid:}
\providecommand{\doi}[1]{\href{http://dx.doi.org/#1}{\path{#1}}}
\providecommand{\Pubmed}[1]{\href{pmid:#1}{\path{#1}}}
\providecommand{\bibinfo}[2]{#2}
\ifx\xfnm\relax \def\xfnm[#1]{\unskip,\space#1}\fi
%Type = Article
\bibitem[{Babuška and Melenk(1997)}]{Babu1997THE}
\bibinfo{author}{Babuška, I.}, \bibinfo{author}{Melenk, J.M.},
  \bibinfo{year}{1997}.
\newblock \bibinfo{title}{The partition of unity method}.
\newblock \bibinfo{journal}{International journal for numerical methods in
  engineering} \bibinfo{volume}{40}, \bibinfo{pages}{727--758}.
%Type = Article
\bibitem[{Babu{\v{s}}ka et~al.(1994)Babu{\v{s}}ka, Caloz and
  Osborn}]{1994Special}
\bibinfo{author}{Babu{\v{s}}ka, I.}, \bibinfo{author}{Caloz, G.},
  \bibinfo{author}{Osborn, J.E.}, \bibinfo{year}{1994}.
\newblock \bibinfo{title}{Special finite element methods for a class of second
  order elliptic problems with rough coefficients}.
\newblock \bibinfo{journal}{SIAM Journal on Numerical Analysis}
  \bibinfo{volume}{31}, \bibinfo{pages}{945--981}.
%Type = Article
\bibitem[{Belytschko and Black(1999)}]{1999Elastic}
\bibinfo{author}{Belytschko, T.}, \bibinfo{author}{Black, T.},
  \bibinfo{year}{1999}.
\newblock \bibinfo{title}{Elastic crack growth in finite elements with minimal
  remeshing}.
\newblock \bibinfo{journal}{International journal for numerical methods in
  engineering} \bibinfo{volume}{45}, \bibinfo{pages}{601--620}.
%Type = Article
\bibitem[{Bingham and Zhang(2007)}]{2007On}
\bibinfo{author}{Bingham, H.B.}, \bibinfo{author}{Zhang, H.},
  \bibinfo{year}{2007}.
\newblock \bibinfo{title}{On the accuracy of finite-difference solutions for
  nonlinear water waves}.
\newblock \bibinfo{journal}{Journal of Engineering Mathematics}
  \bibinfo{volume}{58}, \bibinfo{pages}{211--228}.
%Type = Article
\bibitem[{Chen(2007)}]{chen2007middle}
\bibinfo{author}{Chen, X.B.}, \bibinfo{year}{2007}.
\newblock \bibinfo{title}{Middle-field formulation for the computation of
  wave-drift loads}.
\newblock \bibinfo{journal}{Journal of Engineering Mathematics}
  \bibinfo{volume}{59}, \bibinfo{pages}{61--82}.
%Type = Article
\bibitem[{Chen et~al.(1995)Chen, Molin and Petitjean}]{chen1995numerical}
\bibinfo{author}{Chen, X.B.}, \bibinfo{author}{Molin, B.},
  \bibinfo{author}{Petitjean, F.}, \bibinfo{year}{1995}.
\newblock \bibinfo{title}{Numerical evaluation of the springing loads on
  tension leg platforms}.
\newblock \bibinfo{journal}{Marine Structures} \bibinfo{volume}{8},
  \bibinfo{pages}{501--524}.
%Type = Article
\bibitem[{Cong et~al.(2020)Cong, Teng, Chen and Gou}]{2020A}
\bibinfo{author}{Cong, P.}, \bibinfo{author}{Teng, B.}, \bibinfo{author}{Chen,
  L.}, \bibinfo{author}{Gou, Y.}, \bibinfo{year}{2020}.
\newblock \bibinfo{title}{A novel solution to the second-order wave radiation
  force on an oscillating truncated cylinder based on the application of
  control surfaces}.
\newblock \bibinfo{journal}{Ocean Engineering} \bibinfo{volume}{204},
  \bibinfo{pages}{107278}.
%Type = Inproceedings
\bibitem[{Dai et~al.(2005)Dai, Chen and Duan}]{dai2005computation}
\bibinfo{author}{Dai, Y.S.}, \bibinfo{author}{Chen, X.B.},
  \bibinfo{author}{Duan, W.Y.}, \bibinfo{year}{2005}.
\newblock \bibinfo{title}{Computation of low-frequency loads by the
  middle-field formulation}, in: \bibinfo{booktitle}{20th International
  Workshop for Water Waves and Floating Bodies, Longyearbyen, Norway}, pp.
  \bibinfo{pages}{47--50}.
%Type = Article
\bibitem[{Daux et~al.(2000)Daux, Moes, Dolbow, Sukumar and
  Belytschko}]{2000Arbitrary}
\bibinfo{author}{Daux, C.}, \bibinfo{author}{Moes, N.},
  \bibinfo{author}{Dolbow, J.}, \bibinfo{author}{Sukumar, N.},
  \bibinfo{author}{Belytschko, T.}, \bibinfo{year}{2000}.
\newblock \bibinfo{title}{Arbitrary branched and intersecting cracks with the
  extended finite element method}.
\newblock \bibinfo{journal}{International Journal for Numerical Methods in
  Engineering} \bibinfo{volume}{48}, \bibinfo{pages}{1741--1760}.
%Type = Article
\bibitem[{Engsig-Karup et~al.(2009)Engsig-Karup, Bingham and
  Lindberg}]{engsig2009efficient}
\bibinfo{author}{Engsig-Karup, A.P.}, \bibinfo{author}{Bingham, H.B.},
  \bibinfo{author}{Lindberg, O.}, \bibinfo{year}{2009}.
\newblock \bibinfo{title}{An efficient flexible-order model for {3D} nonlinear
  water waves}.
\newblock \bibinfo{journal}{Journal of Computational Physics}
  \bibinfo{volume}{228}, \bibinfo{pages}{2100--2118}.
%Type = Article
\bibitem[{Espelid and Genz(1994)}]{1994DECUHR}
\bibinfo{author}{Espelid, T.O.}, \bibinfo{author}{Genz, A.},
  \bibinfo{year}{1994}.
\newblock \bibinfo{title}{Decuhr: an algorithm for automatic integration of
  singular functions over a hyperrectangular region}.
\newblock \bibinfo{journal}{Numerical Algorithms} \bibinfo{volume}{8},
  \bibinfo{pages}{201--220}.
%Type = Article
\bibitem[{Fries(2010)}]{2010The}
\bibinfo{author}{Fries, T.P.}, \bibinfo{year}{2010}.
\newblock \bibinfo{title}{{The intrinsic XFEM for two‐fluid flows}}.
\newblock \bibinfo{journal}{International Journal for Numerical Methods in
  Fluids} \bibinfo{volume}{60}, \bibinfo{pages}{437--471}.
%Type = Article
\bibitem[{Fries and Belytschko(2010)}]{2010Fries}
\bibinfo{author}{Fries, T.P.}, \bibinfo{author}{Belytschko, T.},
  \bibinfo{year}{2010}.
\newblock \bibinfo{title}{The extended/generalized finite element method: an
  overview of the method and its applications}.
\newblock \bibinfo{journal}{International journal for numerical methods in
  engineering} \bibinfo{volume}{84}, \bibinfo{pages}{253--304}.
%Type = Article
\bibitem[{Hanssen et~al.(2018)Hanssen, Bardazzi, Lugni and
  Greco}]{hanssen2018free}
\bibinfo{author}{Hanssen, F.C.W.}, \bibinfo{author}{Bardazzi, A.},
  \bibinfo{author}{Lugni, C.}, \bibinfo{author}{Greco, M.},
  \bibinfo{year}{2018}.
\newblock \bibinfo{title}{Free-surface tracking in {2D} with the harmonic
  polynomial cell method: {Two} alternative strategies}.
\newblock \bibinfo{journal}{International Journal for Numerical Methods in
  Engineering} \bibinfo{volume}{113}, \bibinfo{pages}{311--351}.
%Type = Book
\bibitem[{Hughes(2012)}]{hughes2012finite}
\bibinfo{author}{Hughes, T.J.}, \bibinfo{year}{2012}.
\newblock \bibinfo{title}{The finite element method: linear static and dynamic
  finite element analysis}.
\newblock \bibinfo{publisher}{Courier Corporation}.
%Type = Article
\bibitem[{Laborde et~al.(2005)Laborde, Pommier, Renard and
  Sala{\"u}n}]{laborde2005high}
\bibinfo{author}{Laborde, P.}, \bibinfo{author}{Pommier, J.},
  \bibinfo{author}{Renard, Y.}, \bibinfo{author}{Sala{\"u}n, M.},
  \bibinfo{year}{2005}.
\newblock \bibinfo{title}{High-order extended finite element method for cracked
  domains}.
\newblock \bibinfo{journal}{International Journal for Numerical Methods in
  Engineering} \bibinfo{volume}{64}, \bibinfo{pages}{354--381}.
%Type = Inproceedings
\bibitem[{Law et~al.(2020)Law, Liang, Santo, Lim and Chan}]{law2020numerical}
\bibinfo{author}{Law, Y.Z.}, \bibinfo{author}{Liang, H.},
  \bibinfo{author}{Santo, H.}, \bibinfo{author}{Lim, K.Y.},
  \bibinfo{author}{Chan, E.S.}, \bibinfo{year}{2020}.
\newblock \bibinfo{title}{Numerical investigation of the physics of higher
  order effects generated by wave paddles}, in: \bibinfo{booktitle}{Proceeding
  of the 39th International Conference on Ocean, Offshore and Arctic
  Engineering, Fort Lauderdale, FL, USA}, \bibinfo{organization}{American
  Society of Mechanical Engineers Digital Collection}.
%Type = Article
\bibitem[{Li et~al.(2019)Li, Liu and Li}]{2019New}
\bibinfo{author}{Li, A.j.}, \bibinfo{author}{Liu, Y.}, \bibinfo{author}{Li,
  H.j.}, \bibinfo{year}{2019}.
\newblock \bibinfo{title}{New analytical solutions to water wave radiation by
  vertical truncated cylinders through multi-term galerkin method}.
\newblock \bibinfo{journal}{Meccanica} \bibinfo{volume}{54},
  \bibinfo{pages}{429--450}.
%Type = Article
\bibitem[{Liang and Chen(2017)}]{liang2017multi-domain}
\bibinfo{author}{Liang, H.}, \bibinfo{author}{Chen, X.B.},
  \bibinfo{year}{2017}.
\newblock \bibinfo{title}{A new multi-domain method based on an analytical
  control surface for linear and second-order mean drift wave loads on floating
  bodies}.
\newblock \bibinfo{journal}{Journal of Computational Physics}
  \bibinfo{volume}{347}, \bibinfo{pages}{506--532}.
%Type = Article
\bibitem[{Liang et~al.(2015)Liang, Faltinsen and Shao}]{liang2015application}
\bibinfo{author}{Liang, H.}, \bibinfo{author}{Faltinsen, O.M.},
  \bibinfo{author}{Shao, Y.}, \bibinfo{year}{2015}.
\newblock \bibinfo{title}{{Application of a 2D harmonic polynomial cell (HPC)
  method to singular flows and lifting problems}}.
\newblock \bibinfo{journal}{Applied Ocean Research} \bibinfo{volume}{53},
  \bibinfo{pages}{75--90}.
%Type = Article
\bibitem[{Liang et~al.(2020)Liang, Santo, Shao, Law and Chan}]{liang2020liquid}
\bibinfo{author}{Liang, H.}, \bibinfo{author}{Santo, H.},
  \bibinfo{author}{Shao, Y.}, \bibinfo{author}{Law, Y.Z.},
  \bibinfo{author}{Chan, E.S.}, \bibinfo{year}{2020}.
\newblock \bibinfo{title}{Liquid sloshing in an upright circular tank under
  periodic and transient excitations}.
\newblock \bibinfo{journal}{Physical Review Fluids} \bibinfo{volume}{5},
  \bibinfo{pages}{084801}.
%Type = Article
\bibitem[{Ma et~al.(2010a)Ma, Wu and Eatock~Taylor}]{2010FinitePart1}
\bibinfo{author}{Ma, Q.W.}, \bibinfo{author}{Wu, G.X.},
  \bibinfo{author}{Eatock~Taylor, R.}, \bibinfo{year}{2010}a.
\newblock \bibinfo{title}{Finite element simulation of fully non‐linear
  interaction between vertical cylinders and steep waves. part 1: methodology
  and numerical procedure}.
\newblock \bibinfo{journal}{International Journal for Numerical Methods in
  Fluids} \bibinfo{volume}{36}, \bibinfo{pages}{265–285}.
%Type = Article
\bibitem[{Ma et~al.(2010b)Ma, Wu and Eatock~Taylor}]{2010FinitePart2}
\bibinfo{author}{Ma, Q.W.}, \bibinfo{author}{Wu, G.X.},
  \bibinfo{author}{Eatock~Taylor, R.}, \bibinfo{year}{2010}b.
\newblock \bibinfo{title}{Finite element simulations of fully non‐linear
  interaction between vertical cylinders and steep waves. part 2: numerical
  results and validation}.
\newblock \bibinfo{journal}{International Journal for Numerical Methods in
  Fluids} \bibinfo{volume}{36}, \bibinfo{pages}{287--308}.
%Type = Article
\bibitem[{Mavrakos(1988)}]{mavrakos1988}
\bibinfo{author}{Mavrakos, S.}, \bibinfo{year}{1988}.
\newblock \bibinfo{title}{The vertical drift force and pitch moment on
  axisymmetric bodies in regular waves}.
\newblock \bibinfo{journal}{Applied Ocean Research} \bibinfo{volume}{10},
  \bibinfo{pages}{207--218}.
%Type = Article
\bibitem[{Melenk and Babuska(1997)}]{1997Approximation}
\bibinfo{author}{Melenk, J.}, \bibinfo{author}{Babuska, I.},
  \bibinfo{year}{1997}.
\newblock \bibinfo{title}{Approximation with harmonic and generalized harmonic
  polynomials in the partition of unity method}.
\newblock \bibinfo{journal}{Computer Assisted Mechanics and Engineering
  Sciences} \bibinfo{volume}{4}, \bibinfo{pages}{607--632}.
%Type = Article
\bibitem[{Melenk(1995)}]{1995On}
\bibinfo{author}{Melenk, J.M.}, \bibinfo{year}{1995}.
\newblock \bibinfo{title}{On generalized finite element methods}.
\newblock \bibinfo{journal}{PhD thesis, University of Maryland.} .
%Type = Article
\bibitem[{Melenk and Babuška(1996)}]{1996PUFEM}
\bibinfo{author}{Melenk, J.M.}, \bibinfo{author}{Babuška, I.},
  \bibinfo{year}{1996}.
\newblock \bibinfo{title}{The partition of unity finite element method: Basic
  theory and applications}.
\newblock \bibinfo{journal}{Computer Methods in Applied Mechanics and
  Engineering} \bibinfo{volume}{139}, \bibinfo{pages}{289--314}.
%Type = Article
\bibitem[{Moes et~al.(2002)Moes, Gravouil and Belytschko}]{2002Non}
\bibinfo{author}{Moes, N.}, \bibinfo{author}{Gravouil, A.},
  \bibinfo{author}{Belytschko, T.}, \bibinfo{year}{2002}.
\newblock \bibinfo{title}{Non-planar 3d crack growth by the extended finite
  element and level sets. part i : Mechanical model}.
\newblock \bibinfo{journal}{International Journal for Numerical Methods in
  Engineering} \bibinfo{volume}{53}, \bibinfo{pages}{2549--2568}.
%Type = Book
\bibitem[{Newman(2017)}]{Newman2017Marine}
\bibinfo{author}{Newman, J.}, \bibinfo{year}{2017}.
\newblock \bibinfo{title}{Marine Hydrodynamics, 40th Anniversary Edition}.
\newblock \bibinfo{publisher}{MIT Press}.
%Type = Article
\bibitem[{Porter(1995)}]{1995complementary}
\bibinfo{author}{Porter, R.}, \bibinfo{year}{1995}.
\newblock \bibinfo{title}{Complementary methods and bounds in linear water
  waves}.
\newblock \bibinfo{journal}{PhD thesis, University of Bristol.} .
%Type = Book
\bibitem[{Reddy(2019)}]{reddy2019introduction}
\bibinfo{author}{Reddy, J.N.}, \bibinfo{year}{2019}.
\newblock \bibinfo{title}{Introduction to the finite element method}.
\newblock \bibinfo{publisher}{McGraw-Hill Education}.
%Type = Article
\bibitem[{Shao and Faltinsen(2014)}]{shao2014harmonic}
\bibinfo{author}{Shao, Y.}, \bibinfo{author}{Faltinsen, O.M.},
  \bibinfo{year}{2014}.
\newblock \bibinfo{title}{A harmonic polynomial cell (hpc) method for 3d
  laplace equation with application in marine hydrodynamics}.
\newblock \bibinfo{journal}{Journal of Computational Physics}
  \bibinfo{volume}{274}, \bibinfo{pages}{312--332}.
%Type = Inproceedings
\bibitem[{Shao et~al.(2019)Shao, Xiang and Liu}]{shao2019pontoon}
\bibinfo{author}{Shao, Y.}, \bibinfo{author}{Xiang, X.}, \bibinfo{author}{Liu,
  J.}, \bibinfo{year}{2019}.
\newblock \bibinfo{title}{Numerical investigation of wave-frequency pontoon
  responses of a floating bridge based on model test results}, in:
  \bibinfo{booktitle}{Proceeding of the 38th International Conference on Ocean,
  Offshore and Arctic Engineering, Glasgow, Scotland},
  \bibinfo{organization}{American Society of Mechanical Engineers Digital
  Collection}.
%Type = Inproceedings
\bibitem[{Shao et~al.(2016)Shao, You and Glomnes}]{shao2016stochastic}
\bibinfo{author}{Shao, Y.}, \bibinfo{author}{You, J.},
  \bibinfo{author}{Glomnes, E.B.}, \bibinfo{year}{2016}.
\newblock \bibinfo{title}{Stochastic linearization and its application in
  motion analysis of cylindrical floating structure with bilge boxes}, in:
  \bibinfo{booktitle}{Proceeding of the 35th International Conference on Ocean,
  Offshore and Arctic Engineering, Busan, Korea},
  \bibinfo{organization}{American Society of Mechanical Engineers Digital
  Collection}.
%Type = Article
\bibitem[{Shao(2019)}]{2018Numerical}
\bibinfo{author}{Shao, Y.L.}, \bibinfo{year}{2019}.
\newblock \bibinfo{title}{Numerical analysis of second-order mean wave forces
  by a stabilized higher-order boundary element method}.
\newblock \bibinfo{journal}{Journal of Offshore Mechanics and Arctic
  Engineering} \bibinfo{volume}{141}, \bibinfo{pages}{051801}.
%Type = Inproceedings
\bibitem[{Shao and Faltinsen(2012)}]{2012Shao}
\bibinfo{author}{Shao, Y.L.}, \bibinfo{author}{Faltinsen, O.M.},
  \bibinfo{year}{2012}.
\newblock \bibinfo{title}{Towards efficient fully-nonlinear potential-flow
  solvers in marine hydrodynamics}, in: \bibinfo{booktitle}{Proceeding of the
  31st International Conference on Ocean, Offshore and Arctic Engineering},
  \bibinfo{organization}{American Society of Mechanical Engineers}. pp.
  \bibinfo{pages}{369--380}.
%Type = Article
\bibitem[{Shao and Faltinsen(2013)}]{shao2013second}
\bibinfo{author}{Shao, Y.L.}, \bibinfo{author}{Faltinsen, O.M.},
  \bibinfo{year}{2013}.
\newblock \bibinfo{title}{Second-order diffraction and radiation of a floating
  body with small forward speed}.
\newblock \bibinfo{journal}{Journal of offshore mechanics and Arctic
  engineering} \bibinfo{volume}{135}, \bibinfo{pages}{011301}.
%Type = Article
\bibitem[{Strouboulis et~al.(2000a)Strouboulis, Babuška and
  Copps}]{2000Thedesign}
\bibinfo{author}{Strouboulis, T.}, \bibinfo{author}{Babuška, I.},
  \bibinfo{author}{Copps, K.}, \bibinfo{year}{2000}a.
\newblock \bibinfo{title}{The design and analysis of the generalized finite
  element method}.
\newblock \bibinfo{journal}{Computer Methods in Applied Mechanics and
  Engineering} \bibinfo{volume}{181}, \bibinfo{pages}{43--69}.
%Type = Article
\bibitem[{Strouboulis et~al.(2000b)Strouboulis, Copps and
  Babu{\v{s}}ka}]{2000Thegeneralized}
\bibinfo{author}{Strouboulis, T.}, \bibinfo{author}{Copps, K.},
  \bibinfo{author}{Babu{\v{s}}ka, I.}, \bibinfo{year}{2000}b.
\newblock \bibinfo{title}{The generalized finite element method: an example of
  its implementation and illustration of its performance}.
\newblock \bibinfo{journal}{International Journal for Numerical Methods in
  Engineering} \bibinfo{volume}{47}, \bibinfo{pages}{1401--1417}.
%Type = Article
\bibitem[{Sukumar et~al.(2001)Sukumar, Chopp, Moes and
  Belytschko}]{2001Modeling}
\bibinfo{author}{Sukumar, N.}, \bibinfo{author}{Chopp, D.L.},
  \bibinfo{author}{Moes, N.}, \bibinfo{author}{Belytschko, T.},
  \bibinfo{year}{2001}.
\newblock \bibinfo{title}{Modeling holes and inclusions by level sets in the
  extended finite-element method}.
\newblock \bibinfo{journal}{Computer Methods in Applied Mechanics and
  Engineering} \bibinfo{volume}{190}, \bibinfo{pages}{6183--6200}.
%Type = Article
\bibitem[{Sukumar et~al.(2000)Sukumar, Mo{\"e}s, Moran and
  Belytschko}]{2015Extended}
\bibinfo{author}{Sukumar, N.}, \bibinfo{author}{Mo{\"e}s, N.},
  \bibinfo{author}{Moran, B.}, \bibinfo{author}{Belytschko, T.},
  \bibinfo{year}{2000}.
\newblock \bibinfo{title}{Extended finite element method for three-dimensional
  crack modelling}.
\newblock \bibinfo{journal}{International journal for numerical methods in
  engineering} \bibinfo{volume}{48}, \bibinfo{pages}{1549--1570}.
%Type = Article
\bibitem[{Tao et~al.(2007)Tao, Molin and Thiagarajan}]{tao2007spacing}
\bibinfo{author}{Tao, L.}, \bibinfo{author}{Molin, B.and~Scolan, Y.M.},
  \bibinfo{author}{Thiagarajan, K.}, \bibinfo{year}{2007}.
\newblock \bibinfo{title}{Spacing effects on hydrodynamics of heave plates on
  offshore structures}.
\newblock \bibinfo{journal}{Journal of Fluids and structures}
  \bibinfo{volume}{23}, \bibinfo{pages}{1119--1136}.
%Type = Inproceedings
\bibitem[{Taylor and Teng(1993)}]{taylor1993effect}
\bibinfo{author}{Taylor, R.E.}, \bibinfo{author}{Teng, B.},
  \bibinfo{year}{1993}.
\newblock \bibinfo{title}{The effect of corners on diffraction/radiation forces
  and wave drift damping}, in: \bibinfo{booktitle}{Offshore Technology
  Conference, Houston, TX, USA}, \bibinfo{organization}{OnePetro}. pp.
  \bibinfo{pages}{571--581}.
%Type = Article
\bibitem[{Tong et~al.(2021)Tong, Shao, Bingham and Hanssen}]{tong2021adaptive}
\bibinfo{author}{Tong, C.}, \bibinfo{author}{Shao, Y.},
  \bibinfo{author}{Bingham, H.B.}, \bibinfo{author}{Hanssen, F.C.W.},
  \bibinfo{year}{2021}.
\newblock \bibinfo{title}{An adaptive harmonic polynomial cell method with
  immersed boundaries: Accuracy, stability, and applications}.
\newblock \bibinfo{journal}{International Journal for Numerical Methods in
  Engineering} \bibinfo{volume}{122}, \bibinfo{pages}{2945--2980}.
%Type = Article
\bibitem[{Tong et~al.(2019)Tong, Shao, Hanssen, Li and Xie}]{tong2019numerical}
\bibinfo{author}{Tong, C.}, \bibinfo{author}{Shao, Y.},
  \bibinfo{author}{Hanssen, F.C.W.}, \bibinfo{author}{Li, Y.},
  \bibinfo{author}{Xie, B.and~Lin, Z.}, \bibinfo{year}{2019}.
\newblock \bibinfo{title}{Numerical analysis on the generation, propagation and
  interaction of solitary waves by a {Harmonic Polynomial Cell Method}}.
\newblock \bibinfo{journal}{Wave Motion} \bibinfo{volume}{88},
  \bibinfo{pages}{34--56}.
%Type = Article
\bibitem[{Vugts(1968)}]{1968The}
\bibinfo{author}{Vugts, J.H.}, \bibinfo{year}{1968}.
\newblock \bibinfo{title}{The hydrodynamic coefficients for swaying, heaving
  and rolling cylinders in a free surface}.
\newblock \bibinfo{journal}{International Shipbuilding Progress}
  \bibinfo{volume}{15}, \bibinfo{pages}{251--276}.
%Type = Article
\bibitem[{Wu and Eatock~Taylor(1994)}]{wu1994finite}
\bibinfo{author}{Wu, G.}, \bibinfo{author}{Eatock~Taylor, R.},
  \bibinfo{year}{1994}.
\newblock \bibinfo{title}{Finite element analysis of two-dimensional non-linear
  transient water waves}.
\newblock \bibinfo{journal}{Applied Ocean Research} \bibinfo{volume}{16},
  \bibinfo{pages}{363--372}.
%Type = Article
\bibitem[{Wu and Eatock~Taylor(1995)}]{wu1995time}
\bibinfo{author}{Wu, G.X.}, \bibinfo{author}{Eatock~Taylor, R.},
  \bibinfo{year}{1995}.
\newblock \bibinfo{title}{Time stepping solutions of the two-dimensional
  nonlinear wave radiation problem}.
\newblock \bibinfo{journal}{Ocean Engineering} \bibinfo{volume}{22},
  \bibinfo{pages}{785--798}.
%Type = Article
\bibitem[{Xu et~al.(2019)Xu, Zhang, Shao, Gao and Moan}]{xu2019effect}
\bibinfo{author}{Xu, K.}, \bibinfo{author}{Zhang, M.}, \bibinfo{author}{Shao,
  Y.}, \bibinfo{author}{Gao, Z.}, \bibinfo{author}{Moan, T.},
  \bibinfo{year}{2019}.
\newblock \bibinfo{title}{Effect of wave nonlinearity on fatigue damage and
  extreme responses of a semi-submersible floating wind turbine}.
\newblock \bibinfo{journal}{Applied Ocean Research} \bibinfo{volume}{91},
  \bibinfo{pages}{101879}.
%Type = Article
\bibitem[{Yang et~al.(2020)Yang, Teng and Gou}]{2020Comparative}
\bibinfo{author}{Yang, J.}, \bibinfo{author}{Teng, B.}, \bibinfo{author}{Gou,
  Y.}, \bibinfo{year}{2020}.
\newblock \bibinfo{title}{Comparative study on numerical computation methods
  for radiation forces on a three-dimensional body with edge in the time
  domain}.
\newblock \bibinfo{journal}{Journal of Offshore Mechanics and Arctic
  Engineering} \bibinfo{volume}{142}, \bibinfo{pages}{041901}.
%Type = Inproceedings
\bibitem[{Zhao and Faltinsen(1989)}]{zhao1989interaction}
\bibinfo{author}{Zhao, R.}, \bibinfo{author}{Faltinsen, O.M.},
  \bibinfo{year}{1989}.
\newblock \bibinfo{title}{Interaction between current, waves and marine
  structures}, in: \bibinfo{booktitle}{International Conference on Numerical
  Ship Hydrodynamics, 5th, Hiroshima, Japan}, pp. \bibinfo{pages}{513--527}.
%Type = Article
\bibitem[{Zhou and Wu(2015)}]{zhou2015resonance}
\bibinfo{author}{Zhou, B.Z.}, \bibinfo{author}{Wu, G.X.}, \bibinfo{year}{2015}.
\newblock \bibinfo{title}{Resonance of a tension leg platform exited by
  third-harmonic force in nonlinear regular waves}.
\newblock \bibinfo{journal}{Philosophical Transactions of the Royal Society A:
  Mathematical, Physical and Engineering Sciences} \bibinfo{volume}{373},
  \bibinfo{pages}{20140105}.
%Type = Book
\bibitem[{Zienkiewicz et~al.(2005)Zienkiewicz, Taylor and
  Zhu}]{zienkiewicz2005finite}
\bibinfo{author}{Zienkiewicz, O.C.}, \bibinfo{author}{Taylor, R.L.},
  \bibinfo{author}{Zhu, J.Z.}, \bibinfo{year}{2005}.
\newblock \bibinfo{title}{The finite element method: its basis and
  fundamentals}.
\newblock \bibinfo{publisher}{Elsevier}.

\end{thebibliography}

\end{document}